\documentclass[sn-basic]{sn-jnl}

\usepackage{graphicx}%
\usepackage{multirow}%
\usepackage{amsmath,amssymb}
\usepackage{mathrsfs}%
\usepackage[title]{appendix}%
\usepackage{xcolor}%
\usepackage{textcomp}%
\usepackage{manyfoot}%
\usepackage{booktabs}%
\usepackage{algorithm}%
\usepackage{algorithmicx}%
\usepackage{algpseudocode}%
\usepackage{listings}%

\usepackage{epsfig,amsfonts}
\usepackage{amsmath}
\usepackage{amsthm}
\usepackage{graphicx,color}
\usepackage{pgf}
\usepackage{lineno}
\usepackage{mathtools}
\usepackage{adjustbox}
\usepackage{setspace} 
\usepackage{wasysym}
\usepackage{subcaption}

\theoremstyle{thmstyletwo}%
\newtheorem{observation}{Observation}%

\theoremstyle{thmstylethree}%

\newcommand{\cred}[1]{{\color{red}{#1}}}
\newcommand{\cblue}[1]{{\color{blue}{#1}}}

\counterwithin{equation}{section}

\newcommand\Tstrut{\rule{0pt}{2.6ex}}       

\begin{document}

\title[Least squares in low precision arithmetic]{A computational study of low precision incomplete Cholesky factorization preconditioners for sparse linear least-squares problems}


\author*[1,2]{\fnm{Jennifer} \sur{Scott}}\email{jennifer.scott@reading.ac.uk}
\equalcont{These authors contributed equally to this work.}

\author[3]{\fnm{Miroslav} \sur{T\r{u}ma}}\email{mirektuma@karlin.mff.cuni.cz}
\equalcont{These authors contributed equally to this work.}

\affil*[1]{\orgdiv{School of Mathematical, Physical and Computational Sciences}, \orgname{University of Reading}, 
\orgaddress{\city{Reading RG6 6AQ},  \country{UK}}}

\affil[2]{\orgname{STFC Rutherford Appleton Laboratory}, \orgaddress{\city{Didcot, Oxfordshire,\\ OX11 0QX}, \country{UK}.} ORCID iD: 0000-0003-2130-1091}

\affil[3]{\orgdiv{Department of Numerical Mathematics, Faculty of Mathematics and Physics}, 
\orgname{Charles University}, \orgaddress{\country{Czech Republic}.} ORCID iD: 0000-0003-2808-6929}


\abstract{
Our interest lies in the robust and efficient solution of large sparse linear least-squares problems. In recent years, hardware developments have led to a surge in interest in exploiting mixed precision arithmetic within numerical linear algebra algorithms to take advantage of potential savings in memory requirements, runtime and energy use, whilst still achieving the requested accuracy. We explore employing mixed precision when solving least-squares problems, focusing on the practicalities of developing robust approaches using low precision incomplete Cholesky factorization preconditioners. Key penalties associated with lower precision include a loss of reliability and less accuracy in the computed solution. Through experiments involving problems from practical applications, we study computing incomplete Cholesky factorizations of the normal matrix using low precision and using the factors to precondition LSQR using mixed precision. We investigate level-based and memory-limited incomplete factorization preconditioners. We find that the former are not effective for least-squares problems while the latter can provide high-quality preconditioners. In particular, half precision arithmetic can be considered if high accuracy is not required in the solution or the memory for the incomplete factors is very restricted; otherwise, single precision can be used, and double precision accuracy recovered while reducing memory consumption, even for ill-conditioned problems.
}




\maketitle

\section{Introduction}

Depending on the computer architecture, there
can  potentially be significant performance differences when computing and
communicating in different precision formats. Attempts to exploit these differences have resulted in a
long history of efforts to improve the
performance of numerical linear algebra algorithms by seeking to carefully combine
precision formats. The overall goal of mixed precision
algorithms is to accelerate the computational time, and/or to reduce memory requirements
to allow larger problems to be solved or, increasingly, to reduce energy consumption, through 
the judicious employment of
lower precision formats while maintaining robustness and achieving the desired accuracy in the computed solution.
With the growing availability of hardware integration of low precision special function units
that are designed for machine learning applications, classical numerical algorithms are being revisited and
the use of different floating-point formats for performing distinct operations is being explored 
to try and efficiently leverage the available compute resources.
Excellent surveys of 
numerical linear algebra algorithms that seek to exploit mixed precision 
up until 2022 are given in \cite{abdel:2021,hima:2022}; see also the recent discussion in \cite{klbr:2024}.

Our current interest lies in the standard linear least-squares (LS) problem
\begin{equation}\label{eq:ls}
\min_x \|b - A x\|_2,
\end{equation}
where $b \in \mathbb{R}^m$ and $A \in \mathbb{R}^{m \times n}$ are given and we seek $x \in \mathbb{R}^n$.
Our focus is on overdetermined systems ($m > n$). We assume the system matrix $A$ is large, sparse and of full rank but, as is common in practice, it may be ill conditioned.
One solution approach is to employ the normal equations; this may be by forming the normal equations
explicitly or using them implicitly.
It is straightforward to show that $x$ is the unique solution 
of (\ref{eq:ls}) if and only if it satisfies the $n \times n$ {normal equations}
\begin{equation}\label{eq:normal eq}
 C x = A^Tb, \quad C= A^TA . 
\end{equation}
The normal (or Gram) matrix $C$ is symmetric positive definite (SPD) if and only if $rank(A) = n$.
In this case, possible solution approaches
include computing
the Cholesky factorization of $C$ (that is, $C = LL^T$, where $L$
is a lower triangular matrix) or by using a (preconditioned) iterative solver for SPD systems.
However, there are potential problems associated with explicitly forming and using (\ref{eq:normal eq}).
Firstly, information  may be lost when the inner products to compute 
the entries of $C$ are accumulated (see, for example, 
the summary discussion in Chapter 2 of the recent book \cite{bjor:2024}).
Even if the inner products are accumulated in double precision arithmetic, a loss of
information can occur when the computed normal matrix is stored in the working precision.
Indeed, the stored matrix may not be positive definite.
In general, whenever the condition number of $A$ satisfies $\kappa_2(A) \ge u^{-1/2}$
(where $u$ is the machine precision 
and $\kappa_2(A)$ is the ratio of the largest to the smallest singular value of $A$) we can
expect the computed normal matrix to be singular (or indefinite),
in which case computing the Cholesky factorization of $C$ will break down.
Furthermore, although Cholesky factorization algorithms are backward
stable, solution methods that explicitly form the normal equations are not backward
stable because the best backward error bound
contains a factor $\kappa_2(A)$; this is discussed in \cite{high:02}.
Nevertheless, in many practical applications, provided $A$ is not severely ill conditioned, solving the normal equations is regarded as an 
attractive approach, particularly if only modest
accuracy is required; see \cite{hist:87} for a discussion
on why the use of (\ref{eq:normal eq}) can be justified. Note also that in \cite{nart:1992} it is
reported that, for general (square) linear systems, using CG to solve the
corresponding normal equations
is underrated and, despite the squaring of the condition number, it can outperform
other iterative methods applied to the original system.

Iterative methods based on Krylov subspaces are extensively employed
for solving large-scale LS problems. If $A$ is ill conditioned then the
CG (conjugate gradient) method  applied naively to (\ref{eq:normal eq})
can perform poorly. To improve performance,
CGLS \cite{hest:52} uses a slight algebraic rearrangement that avoids explicitly forming $C$. This results in CGLS
having better numerical properties,
with an overhead of some additional  storage and work per iteration.
LSQR (least-squares QR) \cite{pasa:82} is based on the Golub-Kahan (GK) bidiagonalization of $A$. It is also algebraically equivalent 
 to applying CG to (\ref{eq:normal eq}) and (at the cost of extra storage of vectors of length $m$)
is generally more reliable than CGLS when $A$ is ill conditioned and many iterations are needed.

The rate of convergence of Krylov-based  methods depends on the condition number  $\kappa_2(A)$
and on the distribution of the singular values of $A$; convergence 
is generally slow
when $\kappa_2(A)$ is large (many more than $n$ iterations are typically needed to obtain a small backward error), making the use of a preconditioner necessary. 
Krylov solvers can be applied
to the right-preconditioned least-squares problem
\begin{equation*}
\min_{z\in \mathbb{R}^{n}} \|b - AM_R^{-1}z\|_2, \quad x = M_R^{-1}z,
\end{equation*}
by replacing matrix-vector products with $A$ and $A^T$ by products with $AM_R^{-1}$ and
$M_R^{-T}A^T$. The preconditioner must be nonsingular and should be chosen so that $\kappa_2(AM_R^{-1})$ is smaller than $\kappa_2(A)$, or $AM_R^{-1}$ has improved
clustering of its singular values around the origin, and matrix-vector products with $M_R^{-1}$ and $M_R^{-T}$
are relatively inexpensive.
It is well known that finding good preconditioners for LS problems is challenging; see the reviews \cite{bmmt:2014,gosc:2017}.  A common choice is an approximate factorization of
the normal matrix. In particular, if 
we have an incomplete Cholesky (IC)
factorization $C \approx \tilde L \tilde L^T$, where $\tilde L$ is lower triangular, we can set the preconditioner to be
$M_R = \tilde L$.

Our recent studies have shown that, for general SPD systems of equations,
IC factorizations can be computed in low precision and  used to successfully obtain (close to) double precision accuracy
in the final solution, even in the case of highly ill-conditioned systems \cite{sctu:2024,sctu:2025}.
Here we extend this work to an empirical 
investigation of the feasibility of computing and employing low precision incomplete Cholesky (IC) preconditioners for solving sparse least-squares problems. 
The main contributions of this paper are the following.
\begin{itemize}
\item A computational study of stopping criteria for preconditioned LSQR, in particular, implementing and applying the recent work of \cite{pati:2024}.
\vspace{-0.0cm}
\item The robust computation of IC factorizations  of the normal matrix using low precision arithmetic.
\vspace{-0.0cm}
\item A numerical comparison of level-based 
$IC(\ell)$ preconditioners \cite{hypo:02} and memory-limited
preconditioners \cite{sctu:2014a} for LS problems.
\vspace{-0.0cm}
\item An investigation into the effectiveness of low precision IC factorization
preconditioners when used with LSQR and within iterative refinement-based
solvers applied to a range of problems from practical applications, some of which are ill conditioned.
\end{itemize}

We start in Section~\ref{sec:history} with a short overview
of previous work on using mixed precision when solving LS problems. Then, in Section~\ref{sec:test setup}, we describe the test environment to be used in this study
and introduce our set of test problems. In Section~\ref{sec:lsqr},
we recall the iterative solver LSQR and explain the challenging issue of determining when to terminate the
iterations; numerical results are used to support
our choice of stopping criteria.  We also consider incorporating LSQR within a
mixed precision iterative refinement algorithm.
Mixed precision incomplete factorizations of the normal matrix
are presented in Section~\ref{sec:IC factorizations}. The different types of breakdowns that can occur 
during the factorization (particularly when using low precision arithmetic) are discussed, along with the procedures we use to detect potential
breakdowns and then avoid them. This  is key to the development of robust and efficient software that is able to exploit mixed precision arithmetic. Extensive numerical experiments are reported on in Section~\ref{sec:experiments}. These include using a low precision sparse direct solver to compute
a preconditioner as well as incomplete factorization
preconditioners. Finally, in Section~\ref{sec:conclusions}, we summarise our findings and draw some conclusions.

Note that if the sparse matrix $A$ includes one or more rows that are dense (or contain many more
entries than the other rows) then the normal matrix  suffers significant fill-in. This requires modified approaches that
identify and handle such rows separately (see 
the recent papers \cite{sctu:2017b,sctu:2019a,sctu:2022} and references therein). 
Here, we assume all rows of $A$ are sparse
(and only include such problems in our test set).
\newpage
\section{Background}
\subsection{Relation to previous work on solving LS problems using mixed precision}
\label{sec:history}

Previous work on solving linear least-squares problems using mixed precision has been focused on performing complete factorizations in low precision and using them within a variant of iterative refinement.
For square linear systems of equations (both symmetric and non symmetric) there  has been significant interest over the last decade  in the
idea of mixed-precision iterative refinement and, in particular,  GMRES-IR \cite{cahi:17,cahi:18}. The idea is to compute the
matrix factors in low precision and then employ iterative refinement using mixed precision to recover higher precision accuracy.
In GMRES-IR, the correction equation at each refinement step is solved using
the GMRES method preconditioned by the low precision factors.
Although these factors may be of relatively poor quality, they can still
be effective as preconditioners. The analysis shows that 
for matrices that are nearly
numerically singular with respect to the working precision, the condition number of
the preconditioned system is reduced enough to
guarantee backward stability of the approximate solution
computed by preconditioned GMRES.
In contrast, using a basic triangular solve with the low 
precision factors to solve the correction equation
provides no degree of relative accuracy  for even modestly ill conditioned problems. 
Note also that the factors can be used simply as a preconditioner within GMRES, without a refinement loop. The possible disadvantage of this is that the memory requirements and work grows with the number of iterations. The GMRES-IR approach limits the memory and work on each refinement iteration and thus can be viewed as a variant of GMRES with a restarting strategy.
Instead of performing a fixed number of iterations
before each restart, GMRES-IR performs any number of
iterations until it reaches the tolerance set for the GMRES refinement. GMRES without
restarting typically converges faster than with restarting because it uses all of the previously constructed
Krylov subspace to find the new direction to minimize the residual, while the restarted variant starts
constructing a new Krylov subspace after each restart. 

GMRES-IR can be adapted to the least
squares case \cite{hipr:2021}. As least-squares problems,
and the normal equations in particular, may be ill conditioned, iterative
refinement may potentially be used to improve both accuracy and stability. As part of their study of using 
GMRES-IR for symmetric positive definite linear systems, Higham and Pranesh~\cite{hipr:2021}
propose a Cholesky-based GMRES-IR least-squares solver
in which a (complete) Cholesky factorization of the (possibly scaled and shifted) normal equations is computed in low precision and used to compute
an initial approximate least-squares solution. This is then refined to achieve the required accuracy by applying
mixed precision GMRES-IR to the normal equations.
Higham and Pranesh illustrate the potential of
the approach using well conditioned matrices $A$ that are small enough to be handled within their Matlab test environment as dense matrices.

Most recently, Li~\cite{li:2024} investigated the possibility of exploiting mixed precision within the LSQR algorithm when solving discrete linear ill-posed problems
using regularized least-squares. 

An alternative approach to solving least-squares problems using mixed precision  uses the QR factorization of $A$ combined with 
iterative refinement applied to the augmented system formulation \cite{bjor:67a,bjor:2024,cada:2024,cahp:2000,sctu:2022,zhwu:2019}.
This is potentially expensive but more
robust for problems with larger condition numbers.
To reduce the cost, for problems that are highly overdetermined,
modern alternatives to a full QR factorization (such as randomized QR factorizations) have recently been explored \cite{cada:2024a,gbda:2023}. We do not consider using the augmented system in this study.

\subsection{Terminology and test environment}
\label{sec:test setup}
We denote by fp64 and fp32 IEEE double precision (64-bit) and single precision (32-bit), respectively; 
fp16 denotes the 1985 IEEE standard 754 half precision (16-bit).  
Note that bfloat16 is another form of half precision arithmetic but it is not used in our tests
because we use Fortran software (see below) and, as far as we are aware, there are no Fortran compilers that currently support the use of bfloat16. 
Table~\ref{table:precisions} summarises  the parameters for different precision arithmetic. 
We  use $u_{16}$, $u_{32}$, $u_{64}$ to denote the unit roundoffs
in fp16, fp32, and fp64, respectively.
\begin{center}
\begin{table}[htbp]
\small
\begin{tabular}{cllllll}
\hline\Tstrut
&
Signif. &
 Exp. &
$ \;\;u $ &
 $\;x^s_{min}$  &
 $\;x_{min}$  &
 $\;x_{max}$  \\
 \hline\Tstrut
fp16  & 11 & 5 & $4.88 \times 10^{-4}$ & $5.96 \times 10^{-8}$ & $6.10 \times 10^{-5}$ & $6.55 \times 10^4$ \\
fp32 & 24 & 8 & $5.96 \times 10^{-8}$ & $1.40 \times 10^{-45}$ & $1.18 \times 10^{-38}$ & $3.40 \times 10^{38}$ \\
fp64 & 53 & 11 & $1.11 \times 10^{-16}$ &$4.94 \times 10^{-324}$ & $2.22 \times 10^{-308}$ & $1.80 \times 10^{308}$ \\

\hline
\end{tabular}
\caption{Parameters for fp16, fp32, and fp64 arithmetic: the number of bits in the significand and
exponent, unit roundoff $u$, smallest positive (subnormal) number $x^s_{min}$ , smallest normalized positive
number $x_{min}$, and largest finite number $x_{max}$, all given to three significant figures.
}
 \label{table:precisions}
\end{table}
\end{center}

Our test examples are listed in Table~\ref{T:test problems}. They
are  taken from either the  SuiteSparse Matrix Collection\footnote{\url{https://sparse.tamu.edu/}}
or the CUTEst linear
programme set\footnote{\url{https://github.com/ralna/CUTEst}} and comprise a subset of those used by Gould and Scott in their study of numerical methods for
solving large-scale least-squares problems \cite{gosc:2017}. 
Note that some are variants of those in the
SuiteSparse Matrix Collection.
If necessary, the matrix is transposed to give an overdetermined system ($m > n$).
The given condition number estimates  are the ratio of the computed largest and smallest singular values
of $A$. Where available, it is
taken from information provided on the SparseSuite Matrix Collection webpages. Otherwise, an approximation of the largest singular value is computed 
using the iterative procedure proposed in \cite{avdt:2019} (see also \cite{kllu:1996}). 
The smallest singular value is computed using the Matlab routine {\tt svd}. As the condition number of
the normal matrix is the square of the condition number of $A$, we see that the test set
includes problems for which the normal matrix is
highly ill conditioned.

In all our experiments, the matrix $A$ is prescaled so that the 2-norm of each column of 
the scaled matrix $B=AS$ is equal to 1 (here $S$ is the diagonal matrix of scaling factors).
Prescaling is standard within sparse direct solvers (even when using double precision arithmetic) and is also often used with iterative solvers.  
For least-squares problems, scaling corresponds to using a diagonal preconditioner. The study \cite{gosc:2017} demonstrates that this can be highly effective in improving
the performance of an iterative solver. Moreover, it can be beneficially combined with other preconditioners. When using low precision arithmetic,
scaling can reduce the likelihood of overflow and to limit underflows when the matrix is ``squeezed'' into half precision 
(that is, converted from high to low precision) and during the subsequent computation of the preconditioner.
Nevertheless, squeezing the scaled matrix into fp16 can lead to a loss
of information. Higham and Pranesh~\cite{hipr:2021} multiply the scaled matrix by
a scalar so that the resulting normal matrix has constant diagonal entries that do not overflow. For our test matrices with 2-norm scaling,
underflows when converting $B$ from 
fp64 to fp16 only occur for problems co9, psse0 and psse1 and, in each instance, fewer than $1.0 \times 10^2$ entries are lost.

\begin{center}
\begin{table}[htbp]
{
\footnotesize
\begin{tabular}{lrrrrll} 
\hline\Tstrut
{Identifier} &
{$m$} &
{$n$} &
{$nnz(A)$} &
{$nnz(C)$} &
{$density(C)$} &
 {$cond2$} \\ 
\hline\Tstrut
co9  &   22924 &   10789 &$1.10\times 10^5$&$ 1.30\times 10^5$ &$ 2.14\times 10^{-3}$ &$ 5.8\times 10^{4} $ \\
d2q06c &    5831 &    2171 &$3.31\times 10^4$&$  2.92\times 10^4$ &$ 1.19\times 10^{-2}$ &$ 1.4\times 10^{5} $ \\
delf000  &    5543 &    3128 &$1.37\times 10^4$&$ 1.50\times 10^{4}$ &$ 2.74\times 10^{-3}$ &$ 4.2\times 10^{5} $ \\
GE   &   16369 &   10099 &$4.48\times 10^4$&$  6.11\times 10^{4}$ &$ 1.10\times 10^{-3}$ &$ 1.2\times 10^{7} $ \\
IG5-15     &   11369 &    6146 &$3.24\times 10^5$&$  2.87\times 10^{6}$ &$ 1.52\times 10^{-1}$ &$ 2.9\times 10^{19}$ \\
illc1033              &  1033 &  320 &$   4.72\times 10^3$ &$  2.15\times 10^{3}$   &$3.88\times 10^{-2}$ &$  1.9\times 10^4$    \\ 
illc1850       &    1850 &     712 &$8.64\times 10^3$&$ 4.89\times 10^{3}$ &$ 1.79\times 10^{-2}$ &$ 1.4\times 10^{3}$ \\
Kemelmacher &   28452 &    9693 &$1.01\times 10^5$&$  7.24\times 10^{4}$ &$ 1.44\times 10^{-3}$ &$ 2.4\times 10^{4}$  \\
large001  &    7176 &    4162 &$1.89\times 10^4$&$  2.34\times 10^{4}$ &$ 2.46\times 10^{-3}$ &$ 3.9\times 10^{5} $ \\
pilot\_ja  &    2267 &     940 &$1.50\times 10^4$&$  1.53\times 10^{4}$ &$ 3.36\times 10^{-2}$  &$ 2.5\times 10^{8} $ \\
pilotnov  &    2446 &     975 &$1.33\times 10^4$&$  1.31\times 10^{4}$ &$ 2.65\times 10^{-2}$ &$ 3.6\times 10^{9} $ \\
mod2   &   66409 &   34774 &$2.00\times 10^5$&$  3.20\times 10^{5}$ &$ 5.00\times 10^{-4}$ &$ 8.5\times 10^{3} $ \\
psse0      &   26722 &   11028 &$1.02\times 10^5$&$  4.13\times 10^{4}$ &$ 5.88\times 10^{-4}$ &$ 1.0\times 10^{6}$ \\
psse1      &   14318 &   11028 &$5.73\times 10^4$&$  4.51\times 10^{4}$ &$ 6.67\times 10^{-4}$ &$ 2.5\times 10^{8}$ \\
rail2586  &  923269 &    2586 &$8.01\times 10^6$&$ 2.37\times 10^{5}$ &$ 7.05\times 10^{-2}$ &$ 5.0\times 10^{2} $ \\
stat96v2  &  957432 &   29089 &$2.85\times 10^6$&$  1.91\times 10^{5}$ &$ 4.17\times 10^{-4}$ &$ 3.5\times 10^{3} $  \\
watson$\_$1  &  386992 &  201155 &$1.06\times 10^6$&$  1.07\times 10^{6}$ &$ 4.79\times 10^{-5}$ &$ 8.6\times 10^{2} $ \\
well1033              &  1033 &  320 &$   4.73\times 10^3$ &$ 2.15\times 10^{3}$    &$3.88\times 10^{-2}$ &$  1.7\times 10^2$    \\ 
well1850       &    1850 &     712 &$8.76\times 10^3$&$ 4.92\times 10^{3}$ &$ 1.90\times 10^{-2}$ &$ 1.1\times 10^{2}$ \\
world  &   67147 &   34506 &$2.00\times 10^5$&$   3.08\times 10^{5}$ &$ 4.89\times 10^{-4}$ &$ 8.5\times 10^{3} $ \\
\hline
\end{tabular}
}
\caption{Statistics for our test examples. 
$nnz(A)$ and $nnz(C)$ denote the number of entries $A$ and in the lower triangular part of
the normal matrix $C= A^T A$ when computed using fp64 arithmetic. $density(C)$ is the number of nonzeros in $C$ divided by $n^2$.
$cond2$ is an estimate of the 2-norm condition number of $A$.
}
\label{T:test problems}
\end{table}
\end{center}
We take the vector $b$ to be a
vector of random numbers in the interval $[-1,1]$.
This results in the system being inconsistent.
The routine {\tt HSL\_FA14} from the HSL mathematical software library \cite{hsl:2025}
is used to generate $b$.
For a given problem, the same $b$ is used for each experiment.

The software used in our tests is all written in Fortran and compiled using the NAG compiler. As far as we are aware, this is currently the only multi-platform Fortran compiler that fully supports the use of fp16.
The NAG documentation states that their half precision implementation conforms to the IEEE standard. In addition,
using the {\tt -round$\_$hreal} option, all half precision operations are 
rounded to half precision, both at compile time and runtime.
All the conversions this entails result in 
the performance of half precision versions of our software being much slower
than the single precision versions and so reporting time-to-solution is not useful; rather our purpose here is to
simulate low precision arithmetic, enabling us to explore its potential for use in solving tough least-squares problems, without illustrating possible time benefits.

Our experiments involve the iterative solver LSQR (see Section~\ref{sec:lsqr}). We have developed a prototype Fortran implementation that is a modification of the code available 
at \url{ https://web.stanford.edu/group/SOL/software/lsqr/}. An important feature of our version is that it
incorporates a reverse communication user interface that
facilitates the employment of different preconditioners and precisions, and allows
the use of the different stopping criteria discussed in Section~\ref{sec:stopping}. 
It  optionally allows one-sided reorthogonlization; see Section~\ref{sec:intro LSQR}. The reduction in the number of LSQR iterations needed for convergence that results from the use of reorthogonalization is often not sufficient to offset the computational overheads and, as a result, in practical applications LSQR is frequently employed without it. Therefore, unless stated otherwise, reorthogonalization is not used for our reported results.
The HSL package {\tt MI24} is used for experiments with GMRES
and {\tt HSL\_MI35} is used to compute memory-limited IC factorization preconditioners. Note that the software library HSL offers versions of {\tt HSL\_MI35} in fp32 and fp64; for our tests, we developed an fp16 version. For level-based IC($\ell$) factorizations, our implementation is based on the approach of Hysom and Pothen~\cite{hypo:02}). Explicitly forming the normal matrix $C$ is avoided.  We also employ the fp32 and fp64 versions of the sparse direct solver
{\tt HSL\_MA87}~\cite{hors:2010}.

\section{Preconditioned LSQR}
\label{sec:lsqr}
\subsection{Introduction to LSQR}
\label{sec:intro LSQR}

The preconditioned LSQR algorithm \cite{pasa:82} is outlined in Algorithm~\ref{alg:LSQR}. Here $B=AS$ is the scaled LS matrix and the scalars 
$\mu^{(i)} > 0$ and $\beta^{(i)} > 0 $ are chosen to normalise the corresponding vectors; for example,
$\mu^{(1)}p^{(1)} = (BM_R^{-1})^Tq^{(1)}$ implies the computations $\bar p^{(1)} = (BM_R^{-1})^Tq^{(1)}$,
$\mu^{(1)} = \| \bar p^{(1)}\|_2$, $p^{(1)} = (1/\mu^{(1)})\bar p^{(1)}$. In large-scale practical
applications the most expensive part of the computation
is typically the matrix-vector products with
$BM_R^{-1}$ and $(BM_R^{-1})^T$.

\medskip
\begin{algorithm}\caption{Preconditioned LSQR}
\label{alg:LSQR}
\flushleft\textbf{Input:}  $A \in \mathbb{R}^{m \times n}$ and $b\in \mathbb{R}^{m}$,
diagonal scaling matrix $S \in \mathbb{R}^{n \times n} $,
preconditioner $M_R \in \mathbb{R}^{n \times n}$.\\
\flushleft\textbf{Output:} least-squares solution $x$

\setstretch{1.25}\begin{algorithmic}[1]
\State $B = AS$, $z^{(0)} = 0$ \Comment{Scale the matrix}
\State $\beta^{(1)} = \|b\|_2$, $q^{(1)} = b/\beta^{(1)}$
\State $\mu^{(1)}p^{(1)} = (BM_R^{-1})^Tq^{(1)}$,  $w^{(1)} = p^{(1)}$
\State $\bar \rho^{(1)} = \mu^{(1)}$, $\bar \phi^{(1)}= \beta^{(1)}$
\For{$i=1,2, \dots $}
  \State $\beta^{(i+1)}q^{(i+1)} = BM_R^{-1} p^{(i)} - \mu^{(i)}  \,q^{(i)}$\Comment{GK bidiagonalization with preconditioning}
  \State $\mu^{(i+1)}p^{(i+1)}=(BM_R^{-1})^Tq^{(i+1)}- \beta^{(i+1)}  \,p^{(i)}$ 
  \State $\rho^{(i)} = ((\bar \rho^{(i)})^2 + (\beta^{(i+1)})^2)^{1/2}$\Comment{Update QR decomposition}
  \State $c^{(i)} = \bar \rho^{(i)}/\rho^{(i)}$,$\;\; s^{(i)} = \beta^{(i+1)} /\rho^{(i)}$
  \State $\gamma^{(i+1)} = s^{(i)}  \,\mu^{(i+1)}$, $\;\;\bar \rho^{(i+1)} = -c^{(i)}  \, \mu^{(i+1)}$
  \State $\phi^{(i)} = c^{(i)}   \,\bar \phi^{(i)}$ ,
   $\;\;\bar \phi^{(i+1)} = s^{(i)}  \, \bar \rho^{(i)}$
  \State 
 $z^{(i)} \;\;\;\;= z^{(i-1)} +(\phi^{(i)}/\rho^{(i)})  w^{(i)}$ \Comment{Update iterates $z^{(i)}$ and then $w^{(i+1)}$}
  \State $w^{(i+1)} = p^{(i+1)} - (\gamma^{(i+1)}/\rho^{(i)})  w^{(i)}
      $
  \State Test for convergence; exit if converged or maximum iteration count reached
\EndFor
\State $x =  SM_R^{-1}z^{(i)}$ \Comment{Recover LS solution (if not done in Step 14)}
\end{algorithmic}
\end{algorithm}

Steps 6 and 7 of Algorithm~\ref{alg:LSQR} perform the GK bidiagonalization that constructs a Krylov subspace;
Steps 8 to 13 then update the computed solution.
The vectors $p^{(1)}, \ldots ,p^{(i)}$ span an  orthonormal
basis for the Krylov subspace ${\mathcal K}_i(M_R^{-T}B^TBM_R^{-1}, M_R^{-1}Bb) $ and $q^{(1)}, \ldots ,q^{(i)}$ span
${\mathcal K}_i(BM_R^{-1}M_R^{-T}B^T, b) $.
In finite precision arithmetic, the vectors $p^{(i)}$ and $q^{(i)}$ can gradually lose their orthogonality, which can adversely affect convergence, particularly for ill-conditioned problems.  This can be overcome by using
reorthogonalization. This increases the 
work and memory requirements (for large problems, the additional memory needed may be prohibitive, as we illustrate later in  Tables~\ref{T:mi35 fpvar} and \ref{T:mi35 three precisions}) but it can improve convergence
\cite{simo:1984}. When using full reorthogonalization, 
the newly computed vectors $q^{(i+1)}$ and $p^{(i+1)}$ are reorthogonalized against
all previous basis vectors. If $Q^{(i)}$ and $P^{(i)}$
are the matrices of these vectors and are orthonormal to working accuracy, this involves computing
$$
q^{(i+1)} - Q^{(i)} (Q^{(i)})^T q^{(i+1)}, \;\;\;\;
 p^{(i+1)} - P^{(i)} (P^{(i)})^T p^{(i+1)},
$$
and normalizing the resulting vectors (using the Gram Schmidt algorithm). $Q^{(i)}$ and $P^{(i)}$ must be stored, and after $i$ steps the accumulated cost is about $2i^2(m + n)$
flops, making full reorthogonalization impractical for large problems and large $i$. Local (or partial) reorthogonalization limits the costs by reorthogonalizing
$q^{(i+1)}$ and $p^{(i+1)}$
 against a subset of the   previous vectors.
Further savings are made by using one-sided
reorthogonalization in which only the orthonormality of $P^{(i)}$ is maintained \cite{sizh:2000} (see also \cite{fosa:2011,gosc:2017} for numerical results that illustrate the effectiveness of this strategy and \cite{bjor:2024} for further discussion).

\subsection{LSQR stopping criteria}
\label{sec:stopping}
A key issue when developing a practical and robust implementation of an iterative solver is deciding on appropriate stopping criteria. Ideally, we would
like to terminate the computation when
the backward error reaches a user-specified tolerance.
For LS problems, this
may not be straightforward and we have already observed that
LSQR can stagnate. A detailed discussion is given in \cite{cptp:2009}; see also the overview and references
in \cite{hall:2020}.
The linear LS problem we seek to solve can be written as
$$ \min \phi(x), \;\;\; \phi(x) = \|r(x)\|_2 ,\;\;\; r(x) = b - Ax.$$
If the minimum residual is zero ($b \in \mathcal{R}(A)$), $\phi(x)$ is non differentiable at the solution and so
if $x^{(i)}$ is the current  computed solution then 
the first check on its acceptability is on the corresponding residual
$r^{(i)} = b-Ax^{(i)}$. If the minimum residual is nonzero then
$$ \nabla \phi(x) = -\frac{A^Tr(x)}{\|r(x)\|_2}.$$
This leads Gould and Scott~\cite{gosc:2017},
in their comparison study of the performance of different preconditioners for LSQR and LSMR, to use the following stopping rules:
\begin{itemize}
 \item For consistent systems, stop if $\|r^{(i)}\|_2 < \delta_1$.
 \item For inconsistent systems, stop if 
 \begin{equation}\label{eq:ratio_GS}
ratio_{GS} =  \frac{ {\| A^T r^{(i)} \|_2}/{\|r^{(i)}\|_2}} {{\|A^T r^{(0)}\|_2}/{\|r^{(0)}\|_2}} < \delta_2,
 \end{equation}
\end{itemize}
 where $r^{(0)}$ is the initial residual and $\delta_1,\delta_2>0$ are chosen convergence
 tolerances. These criteria for terminating 
 the least-squares solver are  independent of the preconditioner. Thus, they are good for comparing
 preconditioners but may not be appropriate in practice because
 (\ref{eq:ratio_GS}) requires explicitly computing $r^{(i)}$ and $A^Tr^{(i)}$ and thus involves a matrix-vector product with $A$ and $A^T$ each time the computed solution is
 tested for convergence. The overhead can potentially be reduced by  not checking for convergence on every iteration.
 A further issue with using $ratio_{GS}$ as a stopping criteria for LSQR is that, after an initial
 phase in which ${\| A^T r^{(i)} \|_2}/{\|r^{(i)}\|_2}$
 remains constant (or oscillates in ill-conditioned problems), as $i$ increases further this quantity  decreases until
 it and $\|r^{(i)}\|_2$ (and hence also $ratio_{GS}$) stagnate; this is observed and discussed in \cite{cptp:2009}.

In the case of no preconditioning, the above criteria are
closely related to the following from the original LSQR paper~\cite{pasa:82}.
\begin{itemize}
\item For consistent systems, stop if $\|r^{(i)}\|_2 \le \delta_a\,\|A\|_{2,F} \, \|x^{(i)}\|_2 + \delta_b \, \|b\|_2$.
\item For inconsistent systems, stop if
\begin{equation}
    \label{eq:stop paige}
ratio_{PS} = \frac{\|A^{T}r^{(i)}\|_2}{\|A\|_{2,F} \, \|r^{(i)}\|_2} \le \delta_a\,.
\end{equation}
\end{itemize}
Here, $\|A\|_{2,F} $ denotes that either the Frobenius or 2-norm of $A$ may be used. 
The quantities 
$\bar \phi^{(i+1)}$ and $\bar \phi^{(i+1)} \,  \mu^{(i+1)} \, |c^{(i)}|$ within Algorithm~\ref{alg:LSQR} provide estimates of $\|r^{(i)}\|_2 $ and $\|A^Tr^{(i)}\|_2$, which can be used to cheaply check for convergence (no additional
products with $A$ or $A^T$ are needed). 
Furthermore, the $\beta^{(i)}$ and
$\mu^{(i)}$ can be used to accumulate an approximation of $\|A\|_F$; details are given
in \cite{pasa:82}. Observe that when employing
scaling and/or preconditioning,
$A$ is replaced by $BM_R^{-1} = ASM_R^{-1}$ (see Algorithm~\ref{alg:LSQR}) and
using these estimates results in the stopping criteria being based
on  $\|(ASM_R^{-1})^T(b - ASM_R^{-1}z^{(i)})\|_2
= \|(ASM_R^{-1})^T(b - Ax^{(i)})\|_2$,
which depends on $S$ and $M_R$.

The convergence tolerances should ideally be set according to the accuracy of the problem data. If estimates of the relative errors are unknown, 
a small multiple of 
the unit roundoff $u$ may be appropriate \cite{cptp:2009}.
The stopping criteria given above are then sufficient (but not necessary) conditions for $x^{(i)}$ to be
a backward stable LS solution \cite{cptp:2009,jitp:2010}.
However, $\|A^Tr^{(i)}\|_2$ 
can oscillate and it is observed in \cite{cptp:2009} that $\|A^{T}r^{(i)}\|_2 / \|r^{(i)}\|_2$ can plateau and the stopping criteria may not be triggered if 
the tolerances are $ \mathcal{O}(u)$, motivating interest in alternative stopping tests that are applicable with and without preconditioning. Note also that, in practice, the uncertainty in the problem data
and the requirements of the application may mean that much larger
stopping tolerances may be appropriate and more realistic.

From Theorem 5.1 of \cite{cptp:2009},
one possibility is to terminate LSQR using the following condition
\begin{equation}\label{eq:tight stop}
     \frac{\|  x - x^{(i)}\|_{A^TA}}{\|A\|_2 \,\|x^{(i)}\|_2 + \|b\|_2} =\frac{\|P_A \, r^{(i)}\|_2}{\|A\|_2 \, \|x^{(i)}\|_2 + \|b\|_2}  \le \delta,
\end{equation}
where  $x= A^\dag b$ is the solution of the least-squares problem ($A^\dag = (A^TA)^{-1}A^T)$), $\|  x - x^{(i)}\|_{A^TA}$ denotes 
the $A^TA$ norm of the error, that is,  $(x - x^{(i)})^T A^TA(x - x^{(i)})$,
and $P_A = A A^\dag$ is the orthogonal projector onto the range of $A$. It is shown in \cite{cptp:2009} that (\ref{eq:tight stop}) is asymptotically tight in the limit as $x^{(i)}$ approaches the true LS solution and if $\delta = O(u)$ then the computed solution is a backward stable LS solution \cite{cptp:2009,jitp:2010}. 
To use (\ref{eq:tight stop}), estimates of the involved
quantities are needed.
The recent work of Pape\v{z} and Tich\'y~\cite{pati:2024} extends a heuristic-based adaptive estimate proposed originally for the conjugate gradient method \cite{mept:2021} to the iterative solution of the least-squares problem. They
show that the estimate can be evaluated cheaply. Moreover, provided local orthogonality is preserved
during finite-precision computations, and the LSQR solution has not reached its maximal attainable accuracy, the estimate is numerically reliable 
in finite-precision arithmetic, even if the convergence of the
algorithms is strongly influenced by rounding errors.
The estimate is determined dynamically while running LSQR; the adaptive rule used is such that the number of the additional iterations is kept as small as possible.
Importantly, the approach can be used with a split preconditioner and reported results on test examples obtained using Matlab
demonstrate that it is applicable for ill-conditioned problems.

For our experiments, we implement (in Fortran) the adaptive estimate of  \cite{pati:2024}, outlined here as Algorithm~\ref{alg:est}. This takes as input the current iteration $i$, the index $\ell_{i-1}$ determined in the previous iteration and the scalars 
$\{\phi^{(j)}\}_{j=1}^i$ from the LSQR algorithm.
It returns a new index $\ell_i$
and the error norm estimator $estim_{\ell_i}$.
We employ the parameter settings proposed in \cite{pati:2024}, that is,  $\tau = 0.25$
and $tol = 1.0 \times 10^{-4}$.
The latter helps  limit the search for the error estimate to the most significant terms
while $\tau$ represents the relative accuracy of the computed estimate such that
\begin{equation*}
    \frac{\|x - x^{(\ell_i)}\|_{A^TA}^2 - estim_{\ell_i}}{\|x - x^{(\ell_i)}\|_{A^TA}^2} \le \tau .
\end{equation*}
In practice, $\{\Delta_j\}_{j=1}^{i-1}$ can be passed from previous iterations and only $\Delta_i$ is computed in Step 3.   

\begin{algorithm}\caption{Adaptive estimate for LSQR stopping criteria (\ref{eq:tight stop}) at iteration $i$}
\label{alg:est}
\flushleft\textbf{Input:}  Current iteration $i>1$,  $\{\phi^{(j)}\}_{j=1}^i$ from the LSQR algorithm, the index $\ell_{i-1}$ ($\ell_{1}=1$) from iteration $i-1$,  and parameters $\tau$ and $tol$, \ \\
\flushleft\textbf{Output:} $\ell_i $  and error norm estimator  $estim_{\ell_i} \approx  \|x - x^{(\ell_i)}\|_{A^TA}$.

\setstretch{1.17}\begin{algorithmic}[1]
\State $\ell=\ell_{i-1}$
\State $estim_{\ell} = \infty$
\State $\Delta_j=(\phi^{(j)})^2, 1 \le j \le i$ 
\State Set $$p = \text{arg}\max_{j,\, 1 \le j < i}  {\left( \sum_{k=\ell:i} \Delta_k \right)/ \left(\sum_{k=j:i} \Delta_k \right) \le tol};$$ if such $p$ does not exist, set $p=1$
\State Compute $$S = \max_{p \le j < i} \left( \sum_{k=j:i} \Delta_k \right) /\Delta_j$$
\State $\ell_i = \ell$
\While     {$S\, \Delta_i / \left( \sum_{k=\ell:i-1}  \Delta_k\right)  \leq \tau$ and $\ell < i$ }
  \State          $estim_{\ell} = (\sum_{k=\ell:i} \Delta_k)^{1/2}$
  \State          $\ell = \ell + 1$ 
\EndWhile 
\State Set $\ell_i = max(\ell_i,\,\ell-1)$, 
$estim_{\ell_i} = estim_{\ell}$
\end{algorithmic}
\end{algorithm}
\medskip

We propose terminating LSQR when the stopping criteria 
\begin{equation}\label{eq:ratio_PT}
ratio_{PT} = \frac{estim_{\ell_i}}{estim(\|A\|_2) \, \|x^{(i)}\|_2 + \|b\|_2} < \delta
\end{equation}
(computed in fp64 arithmetic) is satisfied for the chosen tolerance $\delta >0$.
Here the estimate $estim(\|A\|_2)$ of the 2-norm of $A$ is computed using the iterative procedure given in \cite{avdt:2019}. In Table~\ref{T:stop test}, we 
present results that compare using the quantities $ratio_{GS}$ (Gould-Scott), $ratio_{PS}$ (Paige-Saunders)
and $ratio_{PT}$ (Pape\v{z}-Tich\'y) given by (\ref{eq:ratio_GS}), (\ref{eq:stop paige})
and (\ref{eq:ratio_PT}), respectively, as the stopping criteria. Here the stopping tolerances (the $\delta$'s) all set to $10^{-10}$. The norm of the final least-squares residual is not reported because each stopping criteria results in effectively the same $\|r\|_2$. The subset of test problems was chosen to illustrate different behaviours. We use the 
single precision variant of the sparse direct linear equation solver {\tt HSL\_MA87} 
\cite{hors:2010} to compute a preconditioner. {\tt HSL\_MA87} is designed to compute
the Cholesky factorization of a sparse symmetric positive definite matrix.
We employ it to compute the Cholesky factorization of the (scaled) normal matrix $AS$. The single precision
$L$ factor is used to precondition LSQR, which is run in double precision (that is, $M_R = L$). 
Each LSQR iteration requires the solution of a 
linear system with $L$ and one with $L^T$. When performing these triangular solves, the single precision factor entries
are locally (in-place) cast to double precision. This requires only a small amount of 
additional double precision memory (of size equal to the largest block on the diagonal of the $L$ factor, see \cite{hors:2010}
for a description of the data structures within {\tt HSL\_MA87}). To accommodate this, 
it was necessary to develop a 
single-double variant of the {\tt HSL\_MA87} solve routine
(currently, there are single and double versions
of {\tt HSL\_MA87} but no mixed precision version
in the HSL Library). We see from Table~\ref{T:stop test} that using the
Gould-Scott test (which is a sufficient but not necessary condition for backward stability), the iteration count
is significantly higher than for the other tests (indeed, LSQR terminates because stagnation occurs without (\ref{eq:ratio_GS}) being satisfied). 
As already observed, the Paige-Saunders test works on the
preconditioned problem. In some cases (including IG5-15 and well1033), 
compared to using the Pape\v{z} and Tich\'y test, this can lead to early termination while for others (such as large001 and psse0) additional 
iterations are performed that result in a smaller $ratio_{PT}$. 


\begin{center}
\begin{table}[htbp]
{\footnotesize
\begin{tabular}{lrrrlrr}
\hline\rule{0pt}{2.6ex}\vspace{-1.5mm} \\ 
\multicolumn{1}{c}{Identifier} & \multicolumn{2}{c}{Gould-Scott} & \multicolumn{2}{c}{Paige-Saunders} & \multicolumn{2}{c}{Pape\v{z}-Tich\'y} \\
\multicolumn{1}{c}{} & \multicolumn{1}{c}{$iters$} &
\multicolumn{1}{c}{$ratio_{PT}$} &
\multicolumn{1}{c}{$iters$} &
\multicolumn{1}{c}{$ratio_{PT}$} &
\multicolumn{1}{c}{$iters$} &
\multicolumn{1}{c}{$ratio_{PT}$} \\
\hline\Tstrut
co9      &  10 &  4.931$\times 10^{-20}$ &  5 &  8.484$\times 10^{-11}$  &  5 &   8.484$\times 10^{-11}$ \\
delf000  &  19 &  1.156$\times 10^{-21}$ & 11 &  4.053$\times 10^{-13}$  &  9 &   5.153$\times 10^{-11}$ \\
IG5-15        &   6 &  9.212$\times 10^{-20}$ & 3  &  3.558$\times 10^{-8}$ &  4 &   2.947$\times 10^{-12}$ \\
large001  &  20 &  2.503$\times 10^{-21}$ & 11 &  2.451$\times 10^{-13}$ &  8 &   7.672$\times 10^{-11}$ \\
mod2      &   8 &  7.202$\times 10^{-22}$ & 4  &  3.799$\times 10^{-11}$ &  4 &   3.799$\times 10^{-11}$ \\
psse0        &  65 &  1.457$\times 10^{-21}$ & 34 &  4.658$\times 10^{-13}$  & 28 &   9.512$\times 10^{-11}$ \\
rail2586  &  14 &  6.612$\times 10^{-18}$ & 8  &  2.308$\times 10^{-10}$ &  9 &   5.706$\times 10^{-11}$ \\
well1033     &   6 &  3.788$\times 10^{-16}$ & 3  &  7.759$\times 10^{-6}$ &   5 &   1.169$\times 10^{-12}$ \\
\hline
\end{tabular}
}
\caption{The effect on the LSQR iteration count of the choice of stopping test.  {\tt HSL\_MA87} run in fp32 arithmetic
is used to compute the preconditioner. For the stopping criteria Gould-Scott (\ref{eq:ratio_GS}), 
Paige-Saunders (\ref{eq:stop paige}), and Papez-Tichy (\ref{eq:ratio_PT}) we report the iteration count $iters$ and $ratio_{PT}$ given by (\ref{eq:ratio_PT})
when LSQR terminates. The stopping tolerances are all set to $10^{-10}$.
}
\label{T:stop test}
\end{table}
\end{center}

\subsection{LSQR-IR} 
As discussed in Section~\ref{sec:history},
Higham and Pranesh \cite{hipr:2021} propose solving linear least-squares problems by applying the GMRES-IR variant of
mixed precision iterative refinement to the normal equations. At each step of iterative refinement, they
employ the low-precision (complete) Cholesky factors of the scaled normal matrix as a preconditioner
for GMRES applied to the correction equation. 
That is, they solve a sequence of linear systems
\begin{equation*}
    A^T A \,d^{\,(i)} = A^T r^{(i)}
\end{equation*}
using preconditioned GMRES. Using GMRES-IR instead of standard iterative refinement enables a much wider range of problems
to be solved \cite{cahi:18}. Obvious variants for
least-squares problems replace  GMRES 
by LSQR and use incomplete Cholesky factorization
preconditioners. The resulting LSQR-IR algorithm using three precisions is given in Algorithm~\ref{alg:lsqr-ir}.
Here, Step 6 is the correction equation.
Note that if $itmax = 1$ and $u_r = u_w$ then the algorithm reduces to two-precision preconditioned LSQR, with the Cholesky factors potentially computed
in a lower precision (for example, we may choose $u_{\ell} = u_{32}$ and $u_w = u_{64}$). If the complete factors
are computed then it is necessary to explicitly forming the normal matrix $B_{\ell}^TB_{\ell}$ (this is generally not necessary for an IC factorization
and is avoided in our numerical experiments). As the factors are used as a preconditioner, the potential loss of information when forming the normal matrix
is less likely to be a concern than would otherwise be the case.

To try and improve efficiency, the use of mixed precision arithmetic within Algorithm~\ref{alg:lsqr-ir}
can be extended by employing mixed precision when solving the
correction equation. 
Following the five-precision variant of GMRES-IR 
proposed in \cite{abhltv:2024}, one possibility is to run LSQR in
precision $u_g \ge u_w$, with the application of the low-precision preconditioner  and products with $A$ and $A^T$ performed in precision $u_p \le u_w$ and the correction $d^{\,(i)}$ stored in precision $u_w$. 
In a recent study of LSQR (without iterative refinement or preconditioning), Li~\cite{li:2024} seeks to achieve potential savings 
by using a two-precision variant of LSQR in which the GK steps and the computation of $z^{(i)}$ and $w^{(i+1)}$ are performed using precision $u_p >u_w$.  The theory given in \cite{li:2024} assumes full reorthogonalization is incorporated within the LSQR algorithm.
If reorthogonalization is needed in practice then this can
add a significant computational overhead (time and memory).

\medskip
\begin{algorithm}\caption{LSQR-IR using three precisions}
\label{alg:lsqr-ir}
\flushleft\textbf{Input:}  $A \in \mathbb{R}^{m \times n}$ and $b\in \mathbb{R}^{m}$,
diagonal scaling matrix $S \in \mathbb{R}^{n \times n} $,
precisions $u_{\ell} \ge u_w \ge  u_r$.\\
\flushleft\textbf{Output:} least-squares solution $x$ in precision $u_w$.
\setstretch{1.17}\begin{algorithmic}[1]
\State Scale $B = AS$ and convert  $B_{\ell} = fl(B)$ in precision $u_{\ell}$
\State  Compute  $C_{\ell}= B_{\ell}^TB_{\ell} \approx LL^T$ in precision $u_{\ell}$\Comment{Either a complete or incomplete factorization}
\State Set $x^{(1)} = 0$ in precision $u_{w}$\Comment{Alternatively, if a complete factorization was computed, use the factors to compute an initial approximate solution $x^{(1)} $}
\For{ $i = 1 : itmax$ or until converged}
\State  Compute the residual $r^{(i)} = b - Ax^{(i)}$ in precision $u_{r}$ and cast it to  $u_{w}$
\State  Solve $\min \|r^{(i)} - Ad^{\,(i)}\|_2$ for the correction $d^{\,(i)}$ using LSQR (Algorithm~\ref{alg:LSQR})
with  preconditioner $M_R = L$ in precision $u_{w}$ 
\State  Update solution $x^{(i+1)} = x^{(i)} + d^{\,(i)}$ in precision $u_{w}$
\State If $x^{(i+1)}$ is sufficiently accurate then set $x=x^{(i+1)}$ and terminate
\EndFor
\end{algorithmic}
\end{algorithm}
\medskip

\subsection{Outer loop termination}

While $ratio_{PT}$ discussed in Section~\ref{sec:stopping} can be used to determine when to terminate each application
of LSQR within LSQR-IR and GMRES-IR, a test 
is needed in the outer refinement loop to decide when to accept the corrected least-squares solution. The backward error
tested in \cite{hipr:2021}
for GMRES-IR applied to the normal equations
is too expensive to compute for large systems. Instead, we  terminate the outer loop when the corrected solution satisfies the stopping criteria (\ref{eq:ratio_GS}) or
$\|r^{(i)}\|_2$ stagnates, that is, 
either $\|r^{(i)}\|_2 > \|r^{(i-1)}\|_2$ or
 \begin{equation}\label{eq:stop outer}
   \frac{ \|r^{(i)}\|_2 - \|r^{(i-1)}\|_2}{\|r^{(i)}\|_2} \le \eta ,
\end{equation}
for a chosen tolerance $\eta \ge 0$.

\section{Incomplete factorization preconditioners in low precision}
\label{sec:IC factorizations}

We start this section by briefly recalling level-based and memory-based incomplete Cholesky (IC) factorizations and then look at when breakdown of the factorization can occur and how to circumvent breakdowns when developing robust implementations, particularly when using low precision arithmetic.

\subsection{Level-based and memory-based IC factorizations}
Incomplete Cholesky factorizations  approximate  the exact Cholesky factorization of a given SPD matrix
$C$ by discarding some entries that occur in a complete factorization. 
Thus $C \approx \widetilde L \widetilde L^T$, where the incomplete
factor $\widetilde L$
is sparse and lower triangular.
The split preconditioned normal equations are 
\[ \widetilde L^{-1} A^T A \widetilde L^{-T} z = \widetilde L^{-1} A^T b, \quad  x = \widetilde L^{-T} z.\]

The simplest sparsity pattern allows no entries
in $\widetilde L$  outside the sparsity pattern 
${\mathcal S}\{C\}$ of $C$.
This 
is termed an $IC(0)$ (or no-fill) factorization.
In practice, sophisticated and systematic ways of extending the sparsity pattern
${\mathcal S}\{\widetilde L\}$ are needed to obtain robust high quality preconditioners.  
One possibility uses the concept of levels \cite{watt:81}. Entries of $\widetilde L$
that correspond to nonzero entries of $C$ are assigned the level 0 while a filled entry in position $(i,j)$ 
(that
is, an entry that is zero in $C$ but nonzero in $L$) is assigned a level as follows:
\begin{equation*}\label{eq:level_definition}
    level(i,j) = \min_{1 \le k < \min\{i,j\}}(level(i,k)+level(k,j)+1).
\end{equation*}
Given $\ell \ge 0$, during the factorization a filled entry  is permitted at position $(i,j)$
provided $level(i,j) \le \ell$.
The number of entries
in $\widetilde L$ (which can be predicted in advance of the numerical factorization using a symbolic factorization) can grow quickly with $\ell$ so only small values  are practical.

Threshold-based incomplete factorizations determine 
the locations of permissible fill-in
in conjunction with the numerical factorization of $C$. Entries of
$\widetilde L$ of absolute value smaller than a prescribed threshold $\tau > 0$ are dropped as they are computed. Unfortunately, choosing a good
$\tau$ can be highly problem dependent.
Memory-based methods prescribe the amount of memory available for $\widetilde L$  and retain only the largest entries in each
row (or column).  Many variations have been proposed (for example, combining
a memory-based method with the use of a drop tolerance); a brief overview is given in \cite{sctu:2011}.

Another possibility is to employ additional memory during the construction of the
incomplete factors that is then discarded. The aim is to obtain a high
quality preconditioner while maintaining sparsity and allowing the
user to control the memory usage \cite{sctu:2014}. 
The Tismenetsky scheme \cite{tism:91} is a matrix decomposition of the form
\begin{equation}\label{eq:memory_limited}
 C = (\widetilde L+ \widetilde R) \, (\widetilde L+ \widetilde R)^T - E,
\end{equation}
where $\widetilde L$ is  lower triangular with positive diagonal
entries, $\widetilde R$ is strictly lower
triangular and $ E = \widetilde R \widetilde R^T$
is the error matrix.
On step $j-1$ ($j \ge 2$) of the factorization, the first column of
the Schur complement is split into the sum
$\widetilde L_{j:n,j}+\widetilde R_{j:n,j},$
where $\widetilde L_{j:n,j}$ 
contains the entries that are retained in  column $j$
of $\widetilde L$\footnote{Here we use the standard section notation, that is, $L_{i:k,j}$ denotes the entries in rows $i$ to $k$ of column $j$}, the diagonal entry $(\widetilde R)_{jj}$ is zero, and $\widetilde R_{j+1:n,j}$
contains the entries that are  not included in $\widetilde L$.
In a complete factorization,
the Schur complement is updated by subtracting
$$(\widetilde L_{j+1:n,j}+\widetilde R_{j+1:n,j}) \, (\widetilde L_{j+1:n,j}+\widetilde R_{j+1:n,j})^T .$$
However, the Tismenetsky incomplete factorization does not compute the full update as it discards the term
\begin{equation*}
E^{(j)} =  \widetilde R_{j+1:n,j} \, \widetilde R_{j+1:n,j}^T .
\end{equation*}
Thus, the matrix $E^{(j)}$ is implicitly added to $C$ and because $E^{(j)}$
is positive semidefinite, in exact arithmetic, the approach does not break down.
An obvious choice is for the largest entries in the column to be retained in $\widetilde L$.
Clearly, more fill entries are used in constructing $L$ than in the standard factorization
and the structure of the complete factorization can be followed more closely. 

Although $\widetilde R$ is discarded once the IC factorization is complete, the columns of $\widetilde R$ must be held until the end of the factorization, independently of the order of operations used by the implementation. 
The computational complexity
and memory needed can be reduced by limiting how many entries are allowed in each column of $\widetilde L$ and $\widetilde R$, as outlined in
Algorithm~\ref{alg:LM_IC}. Drop tolerances may  be optionally used so that only sufficiently large entries are retained in the factors.
This IC factorization algorithm is implemented within the software package {\tt HSL\_MI35}, which is
a modified version of {\tt HSL\_MA28} \cite{sctu:2014,sctu:2014a}. The latter is designed
for general SPD systems while the former is tailored to
solving least-squares problems and optionally avoids explicitly holding the normal matrix.

\medskip
\begin{algorithm}\caption{\bf Left-looking memory-limited IC factorization}\label{alg:LM_IC}
\flushleft\textbf{Input:}  SPD matrix $C\in \mathbb{R}^{n \times n}$ and
${\tt lsize} > 0$ (maximum number of entries in a column of $\widetilde L=\{\tilde l_{ij}\}$) and ${\tt rsize} \ge 0$ (maximum number of entries in a column of $\widetilde R=\{\tilde r_{ij}\}$) \\
\flushleft\textbf{Output:} Incomplete Cholesky factorization  $C \approx  \widetilde L \widetilde L^T$. 
\setstretch{1.17}\begin{algorithmic}[1]
\State Set $w = \{w_i\} = 0$, $1 \le i \le n$\Comment{Initialise work array to zero}
\For{$j=1:n$}\Comment{Loop over the columns}
    \For{$i \in \{ i \ge j \,|\, (i,j) \in \mathcal{S}\{C\}\}$}
     \State    $w_{i} = c_{ij}$ \Comment{Initialise  entries corresponding to nonzeros in $C$}
   \EndFor   
   \For{$k \in \{k < j \, | \, \tilde l_{jk} \neq 0\}$}\Comment{Update column $j$ by column $k$ of $\widetilde L$ if $\tilde l_{jk}\neq 0$}
   \For{$i \in \{ i \ge j \, | \, \tilde l_{ik} \neq 0\}$}
     \State    $w_{i} \leftarrow w_{i}- \tilde l_{ik} \, \tilde l_{jk}$
   \EndFor
   \For{$i \in \{ i \ge j \, | \, \tilde r_{ik} \neq 0\}$}
    \State    $w_{i} \leftarrow w_{i}- \tilde r_{ik} \, \tilde l_{jk}$
    \EndFor
  \EndFor
   \For{$k \in \{k < j \, | \, \tilde r_{jk} \neq 0\}$}\Comment{Update column $j$ by column $k$ of $\widetilde R$ if $\tilde r_{jk}\neq 0$}
   \For{$i \in \{ i \ge j \, | \, \tilde l_{ik} \neq 0\}$}
    \State    $w_{i} \leftarrow w_{i}- \tilde l_{ik}\, \tilde r_{jk}$
        \EndFor
        \EndFor
\State Copy the {\tt lsize} entries of $w$ of largest absolute value into $ \widetilde L_{j:n,j}$
\State Copy  the next largest {\tt rsize} entries of $w$
into $ \widetilde R_{j+1:n,j}$.
\State Scale  $ \tilde l_{jj}=(w_{j})^{1/2}$, \; 
$ \widetilde L_{j+1:n,j} = \widetilde  L_{j+1:n,j}\,/ \tilde l_{jj}\;$,
$ \widetilde R_{j+1:n,j} =  \widetilde R_{j+1:n,j}\,/ \tilde l_{jj}$
\State Reset $w$ to zero. 
  \EndFor
  \State Discard $ \widetilde R$ and return $ \widetilde L$
\end{algorithmic}
\end{algorithm}

\subsection{Avoiding breakdown}

When implementing an IC factorization algorithm 
it is essential to handle the possibility of breakdown.
There are three places where breakdown
can occur, referred to as  B1, B2, and B3 breakdowns   \cite{sctu:2024,sctu:2025}.
\begin{itemize}
\item B1: The diagonal entry $\tilde l_{kk}$ may be unacceptably small or negative. 
\item B2: A column scaling $ \widetilde L_{j+1:n,j} = \widetilde  L_{j+1:n,j}\,/ \tilde l_{jj}\;$,
$ \widetilde R_{j+1:n,j} =  \widetilde R_{j+1:n,j}\,/ \tilde l_{jj}$ may overflow.
\item B3: An update operation $w_{i} \leftarrow w_{i}- \tilde l_{ik} \, \tilde l_{jk}$ or $w_{i} \leftarrow w_{i}- \tilde r_{ik} \, \tilde l_{jk}$ or $w_{i} \leftarrow w_{i}- \tilde l_{ik} \, \tilde r_{jk}$ may overflow.
\end{itemize}
Note that for an $IC(\ell)$ preconditioner, $\widetilde R=0$.
Breakdown can happen when using any precision but is 
most likely for low precision arithmetic.
For higher precision arithmetic, the potential dangers within an
incomplete factorization algorithm can be hidden; 
a standard $IC$ factorization 
using fp64 arithmetic can lead to an ineffective
preconditioner because of
growth in the size of the entries in the factors \cite{sctu:2025}. Without careful monitoring
(which is not routinely done), this growth may be unobserved
but when subsequently applying the preconditioner, very small $\tilde l_{kk}$ can cause the triangular solves to overflow (or come close to overflowing), resulting in the computation aborting or the solver failing to converge. 

Unfortunately, breakdown cannot normally be determined a priori and so the 
development of robust IC factorization implementations 
must seek to avoid breakdowns and to 
detect potential breakdowns as early as possible and then to handle them by revising the data.  
Scott and T\r{u}ma~\cite{sctu:2025} explore a number of
strategies to limit the likelihood of breakdown.
Based on employing $IC(\ell)$ as a preconditioner for a range of SPD problems,
it recommends always prescaling the matrix (which we have already discussed for least-squares problems), incorporating look-ahead, and using a global shifting strategy. We now explain the latter two strategies.

Recall that computing the diagonal entries of
the factor in a (complete or incomplete) Cholesky factorization of an SPD matrix $C=\{c_{ij}\}$ is based on 
\begin{equation*}
    l_{jj} = c_{jj} - \sum_{i < j} l_{ij}^2.
\end{equation*}
Initially, $l_{jj} = c_{jj}$ and at each stage of the factorization a positive (or zero) term
is subtracted from it so that 
$l_{jj}$ either decreases or remains the same.
Thus, to detect potential B1 breakdown as early as possible, look-ahead
can be used whereby, for at each stage $k$ 
the remaining diagonal entries $l_{jj}$ ($j > k$)
are updated and tested. This can be incorporated into a left-looking variant such as is given in Algorithm~\ref{alg:LM_IC}
by holding a separate copy of the diagonal entries.

A consequence of look-ahead 
is that, through the early detection of
B1 breakdowns and taking action to
prevent them,  B3  breakdowns are indirectly prevented.
The numerical experiments in \cite{sctu:2025} on SPD problems coming
from real applications reported that if look-ahead was incorporated then all breakdowns when using
fp16 arithmetic to compute $IC(\ell)$ factorizations were found to be of type B1. Nevertheless,
B2 and B3 breakdowns remain possible and so  the factorization should be implemented using only safe operations,
that is, operations that cannot overflow in the precision being used. 

\begin{figure}[htbp]
\centering
\begin{bordermatrix}{ & k_1 &  k_2 & k_3 &  & j & j_1 &  &  \cr
k_1 &  \ast  \cr
k_2 &     & \ast  &   &  &    \cr
k_3 &  &&   \ast &  \cr
    &   &  &   & \ast &  \cr
j   &  \cblue{\blacksquare} & \cred{\large\CIRCLE} & \cblue{\blacksquare} & & \ast  \cr
j_1  & & \cred{\large\CIRCLE} &  &  &&  \ast \cr
  &   &  &   \cred{\large\CIRCLE} &&e &&  \ast \cr
     & \cblue{\blacksquare} & \cred{\large\CIRCLE} & \cblue{\blacksquare} &&e & && \ast \cr
     & \cblue{\blacksquare} & \cblue{\blacksquare} & \cblue{\blacksquare} &&e& &&& \ast \cr
}
\end{bordermatrix}
        \caption{An example to illustrate the update operations at step $j$ of Algorithm~\ref{alg:LM_IC}.
        Columns $k_1, k_2, k_3$ of $\widetilde L + \widetilde R$ are shown, together with column $j$, which is computed at step $j$.
        Entries denoted by 
squares and circles are computed off-diagonal entries of $\widetilde{L}$
and $\widetilde{R}$, respectively; $e$ denotes an off-diagonal entry of
column $j$ that is updated by one or more of the
columns $k_1, k_2, k_3$.}
        \label{fig:matrix example}
\end{figure}
To demonstrate that, in a left-looking factorization, it
is possible to cheaply check for potential breakdowns, consider the example in Figure~\ref{fig:matrix example}. At step $j$ of Algorithm~\ref{alg:LM_IC}, 
column $j$ of $\widetilde{L}+\widetilde{R}$ is initialised using the entries in column $j$ of $C$ and is then  updated by each column $k < j$ that has a nonzero in row $j$. In our example, $k=k_1, k_2, k_3$;
entries in these columns denoted by 
squares and circles belong to $\widetilde{L}$
and $\widetilde{R}$, respectively. Each update operation
to an entry $(i,j)$ ($i \ge j$)
involves the product of two nonzero entries in column $k$,
one in row $j$ and the other in row $i$. If both 
belong to $\widetilde{R}$ (both are circles)
then the product is excluded. In Figure~\ref{fig:matrix example}, 
the entry $(j_1,j)$ is not updated because the entries $(j,k_2)$ and $(j_1,k_2)$ belong
to $\widetilde{R}$ (circles), whereas the entries denoted
by $e$ in the remaining three rows of column $j$ are updated (if one or more of these is initially zero, then such entries fills in, that is, become nonzero). 
An important feature of the implementation is that 
the entries in the rows of $\widetilde{L}$ and $\widetilde{R}$ are readily available. This means
that the entry of largest absolute value in each row
is also available and can be kept up-to-date
as each column $j$ of $\widetilde L +\widetilde R$ is added, making it straightforward and inexpensive to test for
potential B3 breakdown before the computation of column $j$ commences.

\smallskip
\begin{observation}\label{obs:nobreak}
Let $C$ be a sparse SPD matrix. Assume the
first $j-1$ columns of $\widetilde{L}+\widetilde{R}$  have been successfully computed using precision $u_{\ell}$. 
For  $i \ge j$, let $l_i$ and $r_i$  denote
the number of nonzero entries in $\widetilde{L}_{i,1:j-1}$ and $\widetilde{R}_{i,1:j-1}$, respectively,  and
let $\mu_i$ be  an entry in  of largest absolute in $(\widetilde{L}+\widetilde{R})_{\,i,1:j-1}$.
If $cmax_j$ denotes an entry of largest absolute in column $j$ of $C$ and  $x_{max}$ is the largest finite number in precision $u_{\ell}$, then provided
\begin{equation}\label{eq:assumption}
cmax_j + \mu_i \,\mu_j(\min(l_j,\, \max_{i \ge  j}(l_i+r_i)) + 
\min(r_j,\max_{i \ge  j}l_i) \le x_{max},
\end{equation}
B3 breakdown cannot occur in  step $j$ of Algorithm~\ref{alg:LM_IC}.  
\end{observation}

If (\ref{eq:assumption}) is not satisfied then a more detailed check for B3 breakdown can be performed \cite{sctu:2024}. 
When potential breakdown is detected, we have found that using a global shift is the best approach \cite{sctu:2025}. This
proceeds by choosing a scalar $\alpha > 0$ and
restarting the preconditioner computation by attempting to factorize the scaled and shifted matrix $C(\alpha) = (AS)^T(AS) + \alpha I$
(or, if low precision is used, $C_{\ell}(\alpha) = C_{\ell} + \alpha I$). The hope is that, provided
$\alpha$ is sufficiently small,
the IC factors of $C(\alpha)$ will provide 
an effective preconditioner for the original problem. A simple strategy of repeatedly doubling the shift until the factorization is successful is typically used \cite{hima:2022,limo:99}. However, because an appropriate choice for
the initial shift may not be available, more sophisticated strategies can be beneficial and are used within {\tt HSL\_MI35} \cite{sctu:2014a}.

\section{Numerical experiments}
\label{sec:experiments}

\subsection{Complete factorization preconditioner}
We first consider using a complete factorization
preconditioner, computed by {\tt HSL\_MA87} using single precision.
In Table~\ref{T:lsqr-ir ma87}, we present results for LSQR, LSQR-IR and GMRES-IR using two precisions ($u_r = u_w = u_{64}$, $u_l = u_{32}$).
For LSQR, the stopping criteria is (\ref{eq:ratio_PT}) with $\delta = 10^{-10}$.
For LSQR-IR and GMRES-IR, an initial solution is computed
by solving $LL^T y = SA^Tb $ and setting $x^{(1)}= Sy$
(this is consistent with \cite{hipr:2021}).
For each correction equation, LSQR terminates when
(\ref{eq:ratio_PT}) is satisfied with $\delta = 10^{-5}$;
the convergence tolerance within GMRES is also set to $10^{-5}$. 
The outer iteration is terminated using (\ref{eq:ratio_GS})
with $\delta_2 = 10^{-8}$ and (\ref{eq:stop outer})
with $\eta = 10^3 \times u_{64}$.
We report $ratio_{GS}$ given by (\ref{eq:ratio_GS}), the number $nsol$ of solves with $L$ and $L^T$,
and for LSQR-IR and GMRES-IR the number $nout$ of outer (refinement) iterations.  Note that for LSQR,  $nsol$ is the number of LSQR iterations and this is equal to the number of 
matrix-vector products with $A$ and $A^T$. For LSQR-IR
and GMRES-IR, $nsol$ is one more than
the total number of LSQR and GMRES iterations, respectively, summed over the $nout$ outer iterations (the extra count is from the initial solve), and the number of 
matrix-vector products with $A$ and $A^T$ is $nsol + nout$ (the extra products are required to compute the residual and perform the test (\ref{eq:ratio_GS}) on each outer iteration).
Results are given for a subset of the test set;
the other examples exhibit consistent behaviour.
With the exception of the psse problems, $nsol$ is similar for
all three approaches and, with the high quality
initial solution and complete factorization preconditioner, LSQR-IR and GMRES-IR require only a small number of refinement iterations. Recall that $ratio_{GS}$ provides a sufficient (but not necessary) condition for a backward stable LS solution and, as in Section~\ref{sec:stopping}, the results in Table~\ref{T:lsqr-ir ma87} demonstrate that using
$ratio_{PT}$ as the LSQR stopping criteria can terminate the computation
when $ratio_{GS}$ is  significantly
larger than the tolerance $\delta$. With our parameter choices, the final $ratio_{GS}$
for LSQR-IR and GMRES-IR is typically smaller than for LSQR, suggesting that
the former are potentially over solving the problem and using a smaller
value of the outer loop stopping tolerance $\delta_2$ may be acceptable.

\begin{center}
\begin{table}[htbp]
{
\footnotesize
\begin{tabular}{lrrr rrl rrr}
\hline 
\multicolumn{1}{c}{\Tstrut Identifier} &
\multicolumn{3}{c}{LSQR} &
\multicolumn{3}{c}{LSQR-IR} &
\multicolumn{3}{c}{GMRES-IR} \\
&
\multicolumn{1}{c}{$nsol$} &
\multicolumn{1}{c}{$ratio_{PT}$} &
\multicolumn{1}{c}{$ratio_{GS}$} &
\multicolumn{1}{c}{$nsol$} &
\multicolumn{1}{c}{$nout$} &
\multicolumn{1}{c}{$ratio_{GS}$} &
\multicolumn{1}{c}{$nsol$} &
\multicolumn{1}{c}{$nout$} &
\multicolumn{1}{c}{$ratio_{GS}$} \\
\hline\Tstrut
co9     &    5~~ &   9.010$\times 10^{-11}$ &   2.453$\times 10^{-6}$ &  6~~ & 2~~ &   6.131$\times 10^{-13}$ &    5~~ &     1 &   1.832$\times 10^{-12}$ \\
d2q06c  &    4~~ &   3.181$\times 10^{-11}$ &   8.511$\times 10^{-8}$ &  4~~ & 1~~ &   6.990$\times 10^{-11}$ &    4~~ &     1 &   5.069$\times 10^{-14}$ \\ 
IG5-15  &    4~~ &   3.107$\times 10^{-12}$ &   1.359$\times 10^{-7}$ &  3~~ & 1~~ &   2.442$\times 10^{-12}$ &    3~~ &     1 &   2.210$\times 10^{-12}$ \\  
Kemelmacher & 6~~ &   4.458$\times 10^{-11}$ &   1.892$\times 10^{-7}$ & 6~~ & 2~~ &   9.493$\times 10^{-13}$ &    6~~ &     1 &   5.222$\times 10^{-14}$ \\ 
pilotnov    & 3~~ &   6.369$\times 10^{-12}$ &   5.587$\times 10^{-6}$ & 3~~ & 1~~ &   9.694$\times 10^{-9}$ &    4~~ &     1 &   4.153$\times 10^{-13}$ \\ 
psse0       &28~~ &   9.054$\times 10^{-11}$ &  2.149$\times 10^{-5}$ & 32~~ & 6~~ &   7.767$\times 10^{-9}$ &   44~~ &     1 &   2.289$\times 10^{-10}$ \\
psse1       &25~~ &   3.541$\times 10^{-11}$ &  4.392$\times 10^{-5}$ & 39~~ & 8~~ &   3.881$\times 10^{-9}$ &   30~~ &     1 &   4.135$\times 10^{-10}$ \\ 
rail2586    &10~~ &   6.739$\times 10^{-12}$ &  1.930$\times 10^{-6}$ & 10~~ & 3~~ &   2.923$\times 10^{-10}$ &   13~~ &     2 &   4.317$\times 10^{-15}$ \\
watson\_1   & 4~~ &   4.989$\times 10^{-11}$ &   5.282$\times 10^{-8}$ & 3~~ & 1~~ &   7.901$\times 10^{-12}$ &    3~~ &     1 &   6.664$\times 10^{-12}$ \\  
\hline
\end{tabular}
}
\caption{A comparison of LSQR, LSQR-IR and GMRES-IR. Here {\tt HSL\_MA87} in fp32 arithmetic
is used to compute the Cholesky factorization preconditioner; fp64 is used for the rest of the computation. 
The stopping tolerances are described in the text. $nsol$ denotes the number of solves with $L$ and $L^T$. $ratio_{GS}$ and $ratio_{PT}$
are given by (\ref{eq:ratio_GS}) and (\ref{eq:ratio_PT}), respectively.}
\label{T:lsqr-ir ma87}
\end{table}
\end{center}

In Table~\ref{T:delta ma87}, for LSQR run on problem psse0, we report the ratios
$ratio_{PT}$ and $ratio_{GS}$ for the stopping tolerance $\delta$ in the range $[10^{-6}, \, 10^{-16}]$. The preconditioner is again computed by {\tt HSL\_MA87} in fp32 arithmetic
and applied in fp64 arithmetic.  We see that $ratio_{GS}$ stagnates at $10^{-5}$ while  $ratio_{PT}$
continues to decrease.
Once the maximal attainable accuracy of the LSQR solution has been reached, the scalars $\phi^{(i)}$ are dominated by rounding error rather than by genuine reduction of the error, so the increments $\Delta_i=(\phi^{(i)})^2$ fall to the noise level. Their accumulated sum $estim_{\ell_i}$ (step 8 of Agorithm~\ref{alg:est}), and hence $ratio_{PT}$, can therefore continue to decrease even though $\|x-x^{(i)}\|_{A^TA}$ has stagnated, so the estimate no longer bounds the true error; this is consistent with~\cite{pati:2024}, whose reliability guarantee explicitly excludes this regime.

\begin{center}
\begin{table}[htbp]
{
\footnotesize
\begin{tabular}{lrll| rrll}
\hline 
\multicolumn{1}{c}{\Tstrut $\delta$} &
\multicolumn{1}{c}{$nsol$} &
\multicolumn{1}{c}{$ratio_{PT}$} &
\multicolumn{1}{c|}{$ratio_{GS}$} &
\multicolumn{1}{c}{$\delta$} &
\multicolumn{1}{c}{$nsol$} &
\multicolumn{1}{c}{$ratio_{PT}$} &
\multicolumn{1}{c}{$ratio_{GS}$} \\
\hline\Tstrut
$10^{-6}$   &    5 &   9.572$\times 10^{-7}$ &   1.406$\times 10^{-4}$ &
$10^{-12}$  &   33 &   4.933$\times 10^{-13}$ &   2.441$\times 10^{-5}$ \\ 
$10^{-8}$   &   19 &   6.040$\times 10^{-9}$ &   1.764$\times 10^{-5}$ &
$10^{-14}$  &   40 &   3.787$\times 10^{-15}$ &   1.992$\times 10^{-5}$ \\ 
$10^{-10}$  &   28 &   9.054$\times 10^{-11}$ &   2.149$\times 10^{-5}$ &
$10^{-16}$  &   46 &   3.733$\times 10^{-17}$ &   1.597$\times 10^{-5}$ \\
\hline
\end{tabular}
}
\caption{The effects of varying the stopping tolerance $\delta$ on $ratio_{PT}$ and $ratio_{GS}$
given by (\ref{eq:ratio_PT}) and (\ref{eq:ratio_GS}), respectively. LSQR is preconditioned by the
Cholesky factor $L$ computed using {\tt HSL\_MA87} in fp32 arithmetic. The test problem is psse0. $nsol$ denotes the number of solves with $L$ and $L^T$.}
\label{T:delta ma87}
\end{table}
\end{center}

\subsection{Level-based incomplete factorization preconditioner}

\begin{figure}
     \centering
     \begin{subfigure}[b]{0.4\textwidth}
         \centering
         \includegraphics[width=\textwidth]{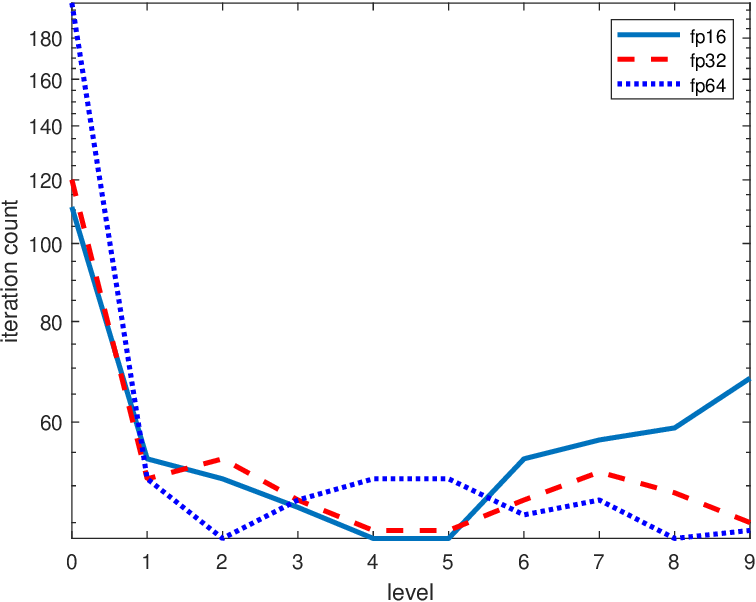}
         \caption{Problem GE}
     \end{subfigure}
     \hspace{1cm}
     \begin{subfigure}[b]{0.4\textwidth}
         \centering
         \includegraphics[width=\textwidth]{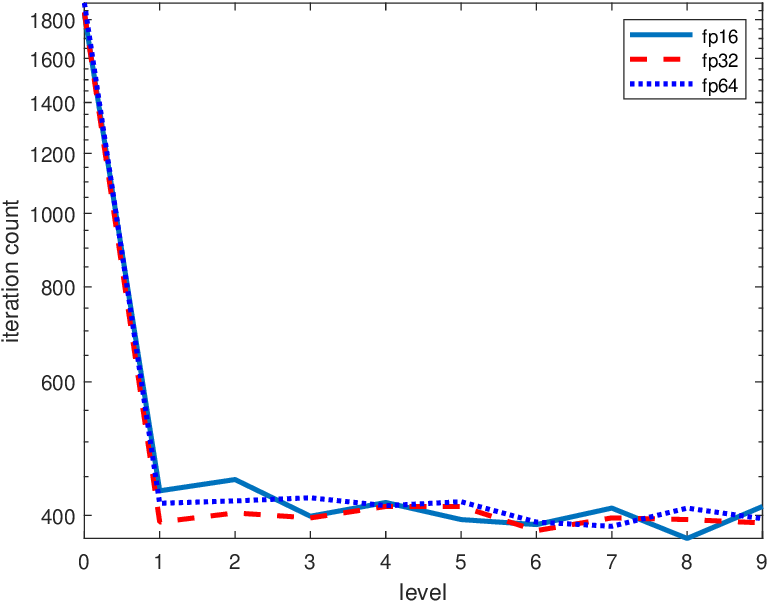}
         \caption{Problem mod2}
     \end{subfigure}
        \caption{LSQR iteration counts for problems GE (left)
        and mod2 (right) as the level parameter $\ell$  increases. The IC factorization is computed fp16, fp32 and fp64 arithmetic. The stopping criteria  is (\ref{eq:ratio_PT}) with $\delta = 10^{-5}$.}
        \label{fig:IC(l)}
\end{figure}

In \cite{sctu:2024,sctu:2025},  we explored avoiding breakdown when computing level-based $IC(\ell)$ preconditioners
in low precision arithmetic. For general sparse SPD linear systems from a variety of practical applications
we found that, when carefully implemented, it was possible to compute the preconditioner using half precision arithmetic and to employ it within a Krylov subspace-based refinement algorithm (for example, GMRES-IR)
to recover double precision accuracy in the computed solution. Moreover, the increase in iteration count resulting from the use of fp16 was generally only significant for highly ill-conditioned examples.
Thus, it is of interest to consider whether 
$IC(\ell)$ preconditioners can be effective for LS problems.
Representative results are presented  in Table~\ref{T:IC(l)} for LSQR preconditioned by the $IC(\ell)$ factor with $\ell = 3$ computed using fp16 and fp64 arithmetic.
Here and elsewhere, fp1/fp2 denotes the incomplete factorization is computed using fp1 arithmetic and fp2 arithmetic is used for matrix-vector products with $A$ and $A^T$ and for applications of the preconditioner.
We see that, even if fp64 arithmetic is used throughout,
in many cases the $IC(3)$ preconditioner performs poorly and
LSQR fails to converge for small values of the stopping tolerance $\delta$; similar disappointing convergence is seen  for the other problems in our test set and so the results for them are omitted from the table.
Experiments with other choices of $\ell \le 5$ also fail to result
in effective preconditioners for use with LSQR. This is illustrated
in Figure~\ref{fig:IC(l)}. Note in particular that $IC(0)$ (which allows only entries corresponding to the entries in $A^TA$) requires
a large number of iterations to obtain the requested accuracy.

\begin{center}
\begin{table}[htbp]
{
\small
\begin{tabular}{lr rrr| rrr}
\hline\Tstrut
& &
\multicolumn{3}{c|}{fp16/fp64} &
\multicolumn{3}{c}{fp64/fp64} \\

\multicolumn{1}{c}{\Tstrut Identifier} &
\multicolumn{1}{c}{$nz(\tilde L)$} &
\multicolumn{1}{c}{$10^{-5}$} &
\multicolumn{1}{c}{$10^{-7}$} &
\multicolumn{1}{c|}{$10^{-9}$}  &
\multicolumn{1}{c}{$10^{-5}$} &
\multicolumn{1}{c}{$10^{-7}$} &
\multicolumn{1}{c}{$10^{-9}$}  \\
\hline\Tstrut
co9       & $4.98 \times 10^5$ & 2179 & $\dag$ & $\dag$ & 2195 & $\dag$ & $\dag$ \\
delf000   & $4.18 \times 10^4$ & 372  & $\dag$ & $\dag$ & 1675 & $\dag$ & $\dag$ \\
GE        & $1.99 \times 10^5$ &  51  & $\dag$ & $\dag$ &  46 & $\dag$ & $\dag$ \\
IG5-15    & $1.20 \times 10^7$ & 547  & 1229  & 1720  &  623 & 1247 &  1689 \\
large001  & $7.15 \times 10^4$ & 767  & $\dag$ & $\dag$ & 753 & $\dag$ & $\dag$ \\
mod2      & $1.86 \times 10^6$ & 432  & 2174 & $\dag$ &  417   & 2095 &  $\dag$  \\
pilot\_ja  & $6.27 \times 10^4$ &  4   &  11  & $\dag$ &  4    & 13 & $\dag$ \\
psse0     & $5.58 \times 10^4$ & 401  & $\dag$ & $\dag$ & $\dag$ & $\dag$ & $\dag$ \\ 
rail2586  & $1.15 \times 10^6$ & 568  & 1004 & 1123 &  531  & 901   & 1286 \\
well1033  & $2.77 \times 10^3$ & 312  & 339  & 357 &  315 &  333 &   360\\
\hline
\end{tabular}
}
\caption{Iteration counts for LSQR preconditioned by the
$IC(3)$ factor computed using fp16 and fp64 arithmetic
for a range of stopping tolerances.
$nz(\tilde L)$ denotes the number of
entries in the incomplete factor. The stopping criteria (\ref{eq:ratio_PT}) is used with $\delta= 10^{-5},10^{-7}, 10^{-9}$. 
$\dag$ indicates the stopping criteria is not achieved within 3000 LSQR iterations.
}
\label{T:IC(l)}
\centering
\end{table}
\end{center}

\subsection{Memory-limited incomplete factorization preconditioner}
We next study the behaviour of memory-limited IC preconditioners.
Table~\ref{T:mi35 fpvar} presents results for the IC preconditioner
computed using {\tt HSL\_MI35}. 
Here the parameter {\tt lsize} is chosen so that the
number $nz(\tilde L)$ of entries in the incomplete factor is similar to 
that for the level-based $IC(3)$ preconditioner
that was reported on in the previous section (with a maximum value of 60). The parameter {\tt rsize} controls the
number of entries in the temporary factor $\widetilde R$ 
used in the construction of $\widetilde L$  (recall Algorithm~\ref{alg:LM_IC}). In all our tests, we set
${\tt rsize} = {\tt lsize}$. A number of observations can be made.
Firstly, it clear that the memory-limited preconditioner
is much more robust compared to the $IC(\ell)$ preconditioner.
With the stopping tolerance of $10^{-5}$, for many of our
test examples, the iteration counts when employing the fp16 preconditioner
are competitive with those for the fp32 and fp64 preconditioners
(exceptions include illc1033 and rail2586).
If greater accuracy is requested ($\delta = 10^{-10}$) then,
with the chosen {\tt lsize} settings,
the fp16 preconditioner requires 
significantly more iterations
than the higher precision preconditioners for many (but not all) examples. Incorporating reorthogonalization 
can reduce the iteration counts and, in some instances (including illc1033, psse0 and psse1), the reduction is substantial, although
the counts typically remain high (see, for instance, rail2586). Further investigation of the performance of the fp16 preconditioner for the psse problems
reveals that, when compared to using $u_{64}$, a significant number of entries (approximately $7.8 \times 10^3$)
are lost when each column of the normal matrix $C_l = B_l^TB_l$ is computed in precision $u_l  = u_{16}$ ({\tt HSL\_MI35} allows the columns to be computed as required, rather than computing and storing $C_l$). 
This can be regarded as an initial sparsification of the psse normal matrix, which
may account for the poor quality preconditioner.
For all the other test problems,
either no entries or a very small number of entries are lost. To try and reduce the loss of information for the psse problems,
we experimented with the strategy of Higham and Pranesh~\cite{hipr:2021}, which 
scales all the entries of $B = AS$ by a factor $\mu^{1/2}$, where $\mu =  \theta x_{max}$ and $\theta <1$ is a chosen parameter, before
the matrix is squeezed into fp16. Setting $\hat C_l= (\mu^{1/2} AS)^T(\mu^{1/2} AS)$,
Higham and Pranesh also propose using the scaled and shifted normal matrix 
$\hat C_l + c\,u_l\,D$, where $D$ is a diagonal matrix with entries $(\hat C_l)_{ii}$
and $c>0$ is a chosen constant. In exact arithmetic, the diagonal entries of $\hat C_l$ are all equal to
$\mu$ and overflow cannot happen for $\theta < 1$. In our
experiments, we found that scaling $AS$ with $\mu^{1/2}$ 
did not substantially reduce the
lost entries when squeezing the matrix into fp16.  For some test cases, the LSQR iteration count was reduced by using this approach but for other examples, using the simpler shifted matrix $C_l + \alpha I$ gave a lower count.
As the latter requires fewer parameters to be chosen, we have used it in all our experiments (which is consistent with the fp32 and fp64 experiments).

\begin{center}
\begin{table}[htbp]
{
\footnotesize
\begin{tabular}{lrr rrr| rlrr}
\hline\Tstrut
 & & &
\multicolumn{3}{c|}{$\delta = 10^{-5}$} &
\multicolumn{3}{c}{$\delta = 10^{-10}$}  \\

\multicolumn{1}{c}{\Tstrut Identifier} &
\multicolumn{1}{c}{{\tt lsize}} &
\multicolumn{1}{c}{$nz(\tilde L)$} &
\multicolumn{1}{c}{fp16} &
\multicolumn{1}{c}{fp32} &
\multicolumn{1}{c|}{fp64} &
\multicolumn{2}{c}{fp16} &
\multicolumn{1}{c}{fp32} &
\multicolumn{1}{c}{fp64} \\
\hline\Tstrut
co9       & 45 & $4.95 \times 10^5$ &  38 & 15 & 12 & 277 & (224) & 93 & 119 \\
d2q06c    & 60 & $1.29\times 10^5$ &   8  &  3 &   4 &  71 &(71) & 10 & 10 \\
delf000   & 20 & $6.53\times 10^4$ &   14  &  4 &   3 & 658 & (622) & 17 & 5 \\
GE        & 15 & $1.61 \times 10^5$ &   4 & 2 & 5 & 219 &(203) & 27 & 40 \\
IG5-15    & 60 & $ 3.73 \times 10^5$ &   98 & 98 & 89 & 305 &(290) & 305 & 301 \\
illc1033  &  10 & $2.97 \times 10^3$ &  245 & 15 &   3 & 305 & (156) & 19 & 3 \\
illc1850  &  10 & $7.76 \times 10^3$ &   82 & 30 &  32 & 126 & (94) & 36 & 39 \\
Kemelmacher& 30 & $ 2.99 \times 10^5$ &  46 &  43 &  15 & 72 & (72) & 68 & 66 \\
large001  & 15 & $ 6.64 \times 10^4$ &   26 &  4  &   4 & 835 & (805) & 17 & 17 \\
mod2      & 30 & $1.08 \times 10^6$&   10  &  7  &  5 & 89 & (81) & 89 & 84 \\
pilot\_ja  & 60 & $ 5.51 \times 10^4$ &   2  &  2  &   2 & 42 & (42) & 5 & 5 \\
pilotnov  & 60 & $ 5.69 \times 10^4$ &   2  &  2  &   2 & 22 & (22) & 3 & 5 \\
psse0     & 5 & $6.61 \times 10^4$&    16 &  5  &   2 & $\dag $ & (2174) & 107 & 47 \\ 
psse1     & 60 & $ 5.09 \times 10^5$ &   8  &  5  &   2 & $\dag$ & (1829) & 64 & 3 \\ 
rail2586  & 60 & $ 1.53 \times 10^5$ &  333 & 80 &  60 & 800 & (642) & 111 & 103 \\
stat96v2  & 20 & $ 6.11 \times 10^5$  & 37 & 10  & 10 &  64 & (64) & 18 & 18 \\
watson\_1 & 15 & $ 3.22 \times 10^6$  & 59 & 52  & 53 & 163 & $~\ast$ & 141 & 142 \\
well1033  & 10 & $ 2.95 \times 10^3$   & 12 &  3  &  3 &  17 & (17) & 5 & 3 \\
well1850  & 10 & $ 7.77 \times 10^3$   & 11 & 11  & 12 &  19 & (19) & 18 & 19 \\
world     & 50 &  $ 1.76 \times 10^6$ & 8  &  3  & 3 &  62 & (61) & 17 & 16 \\
\hline
\end{tabular}
}
\caption{Iteration counts for LSQR preconditioned by the
memory-limited IC factor computed  by {\tt HSL\_MI35} using
fp16, fp32 and fp64 arithmetic. All matrix-vector products with $A$ and applications of the preconditioner are performed using fp64 arithmetic. $nz(\tilde L)$ denotes the number of
entries in the incomplete factor. The stopping criteria is (\ref{eq:ratio_PT}) with $\delta= 10^{-5}$ and $10^{-10}$. 
$\dag$ indicates requested accuracy not achieved within 3000 LSQR iterations. The numbers in parentheses are iteration
counts for LSQR with full 
one-sided reorthogonalization. $\ast$ indicates 
insufficient memory. 
}
\label{T:mi35 fpvar}\vspace{3mm}
\end{table}
\end{center}

Finally, we observe that if fp16 (respectively, fp32) is used
in place of fp64, the memory saving for each entry of $\tilde L$ is 6 (respectively, 4) bytes. Thus, for problems for which $nz(\tilde L)$ is
large compared to $nz(A)$ (problems that are not tall and skinny), significant memory savings result from using low precision.

Table~\ref{T:delta mi35} illustrates varying the stopping tolerance $\delta$.
Here, for the moderately ill-conditioned test problem psse0 (its estimated condition number is $10^6$) the IC factorization preconditioner is computed by {\tt HSL\_MI35} using fp32 arithmetic. 
We see that
$ratio_{GS}$ stagnates at around $10^{-12}$ while $ratio_{PT}$ decreases further; similar findings are observed for
other test examples. Recall that $ratio_{PT}$ is only reliable until the 
level of maximal
attainable accuracy is reached, but in practice this level is not known.
\begin{center}
\begin{table}[htbp]
{
\footnotesize
\begin{tabular}{lrll| rrll}
\hline
\multicolumn{1}{c}{\Tstrut $\delta$} &
\multicolumn{1}{c}{$nsol$} &
\multicolumn{1}{c}{$ratio_{PT}$} &
\multicolumn{1}{c|}{$ratio_{GS}$} &
\multicolumn{1}{c}{$\delta$} &
\multicolumn{1}{c}{$nsol$} &
\multicolumn{1}{c}{$ratio_{PT}$} &
\multicolumn{1}{c}{$ratio_{GS}$} \\
\hline\Tstrut
$10^{-6}$  &    13 &   9.523$\times 10^{-7}$ &  2.229$\times 10^{-4}$ &
$10^{-12}$  &   128 &   7.837$\times 10^{-13}$ &  9.997$\times 10^{-12}$ \\ 
$10^{-8}$  &   90 &   6.040$\times 10^{-9}$ &   1.764$\times 10^{-5}$ &
$10^{-14}$   &   147 &   9.532$\times 10^{-15}$ &  8.212$\times 10^{-13}$ \\ 
$10^{-10}$  &   107 &   9.446$\times 10^{-11}$ &  1.502$\times 10^{-9}$ &
$10^{-16}$   &   164 &   9.755$\times 10^{-17}$ &  4.824$\times 10^{-13}$ \\

\hline
\end{tabular}
}
\caption{The effects of varying the stopping tolerance $\delta$ on $ratio_{PT}$ and $ratio_{GS}$
given by (\ref{eq:ratio_PT}) and (\ref{eq:ratio_GS}), respectively. LSQR is preconditioned by the
IC factor $\tilde L$ computed using {\tt HSL\_MI35} in fp32 arithmetic. The test problem is psse0 with ${\tt lsize}=5$. $nsol$ denotes the number of solves with $L$ and $L^T$.}
\label{T:delta mi35}\vspace{3mm}
\end{table}
\end{center}

\begin{figure}
     \centering

        \begin{subfigure}[b]{0.4\textwidth}
         \centering
         \includegraphics[width=\textwidth]{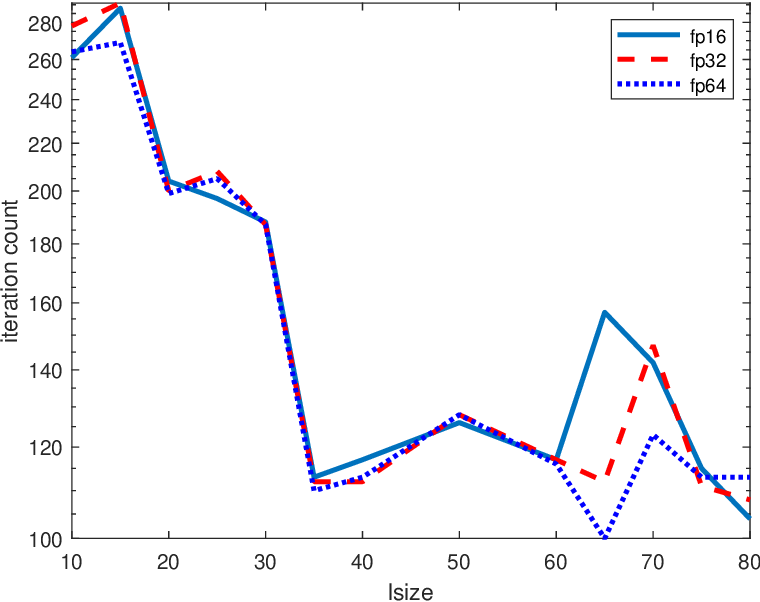}
         \caption{Problem IG5\_15 with $\delta = 10^{-5}$}
     \end{subfigure}
     \hspace{1cm}
     \begin{subfigure}[b]{0.4\textwidth}
         \centering
         \includegraphics[width=\textwidth]{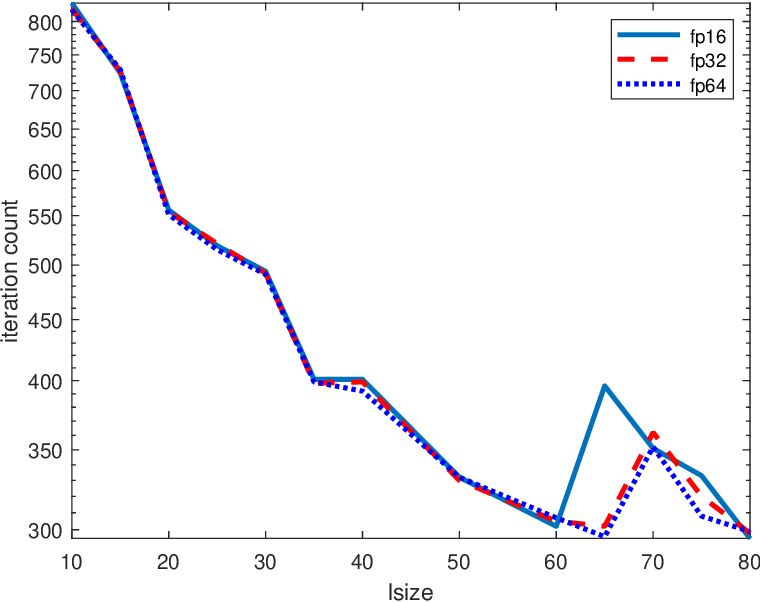}
         \caption{Problem IG5\_15 with $\delta = 10^{-10}$}
     \end{subfigure}

     \vspace{1cm}
     \begin{subfigure}[b]{0.4\textwidth}
         \centering
         \includegraphics[width=\textwidth]{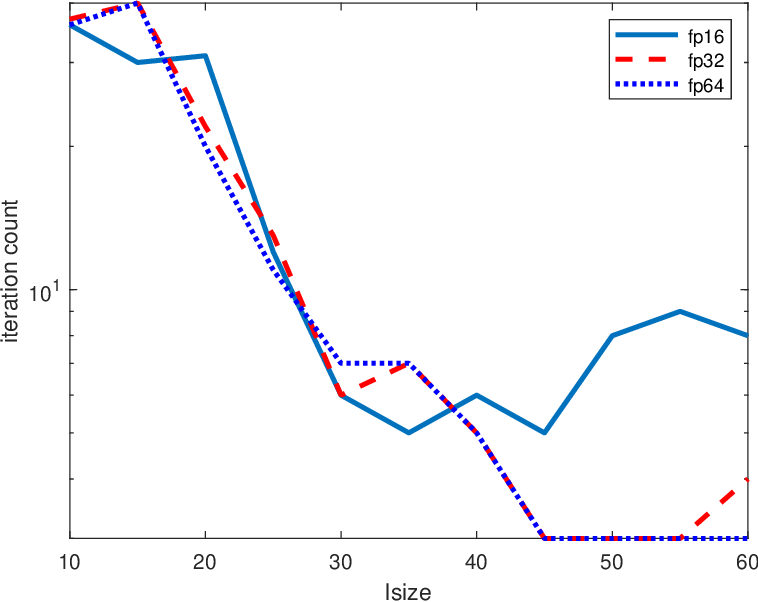}
         \caption{Problem world with $\delta = 10^{-5}$}
     \end{subfigure}
     \hspace{1cm}
     \begin{subfigure}[b]{0.4\textwidth}
         \centering
         \includegraphics[width=\textwidth]{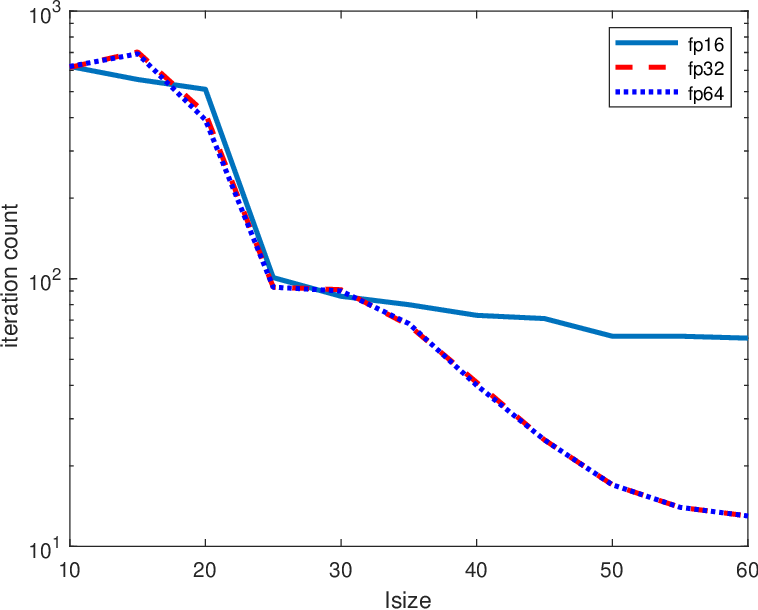}
         \caption{Problem world with $\delta = 10^{-10}$}
     \end{subfigure}

        \caption{LSQR iteration counts for problems IG5\_15 (top) and world (bottom) as the {\tt HSL\_MI35} parameter {\tt lsize} that controls the number of entries in the IC factor is increased. The IC factorization is computed fp16, fp32 and fp64 arithemtic. The stopping criteria is (\ref{eq:ratio_PT}).}
        \label{fig:vary lsize}
\end{figure}

An important feature of the memory-limited IC factorization is that,
via the parameter {\tt lsize},
the user can control the memory used and can limit the number
of entries in each column of $\tilde L$.
In many instances, as {\tt lsize} increases, so does the quality of
the resulting preconditioner (but the larger number of entries in
$\tilde L$ not only requires more memory but also increases the cost of each application of the preconditioner). In Figure~\ref{fig:vary lsize}, 
for problems IG5\_15 and world, we plot the
LSQR iteration count as {\tt lsize} increases.

\begin{figure}
     \centering
     \begin{subfigure}[b]{0.4\textwidth}
         \centering
         \includegraphics[width=\textwidth]{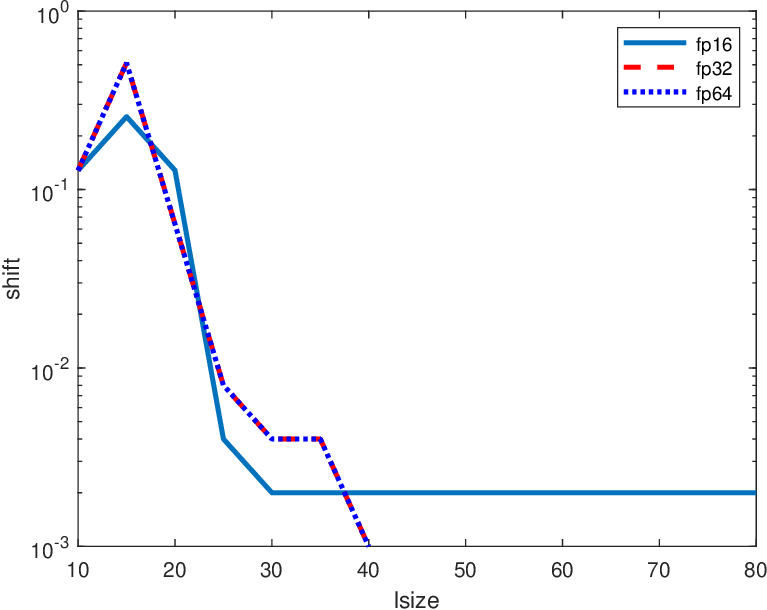}
        \caption{The shift $\alpha$ as {\tt lsize} increases}
     \end{subfigure}

     \vspace{1cm}

          \begin{subfigure}[b]{0.4\textwidth}
         \centering
         \includegraphics[width=\textwidth]{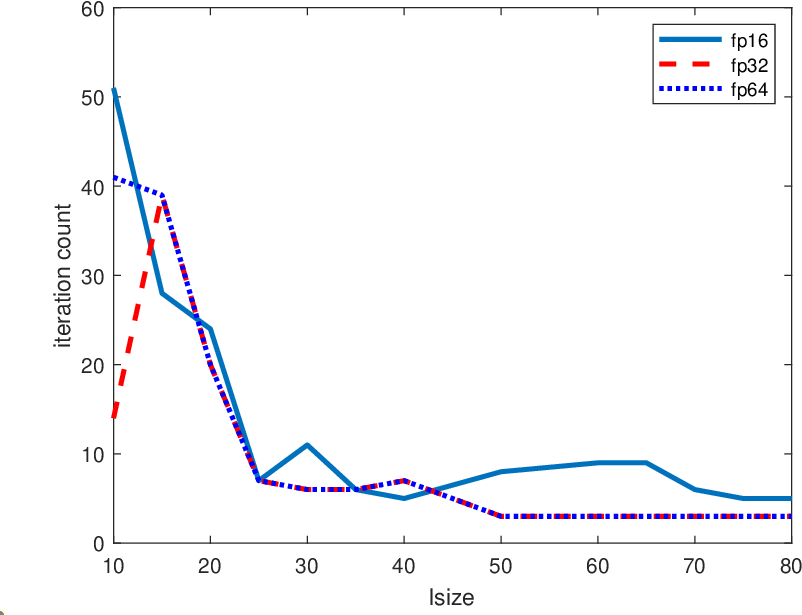}
       \caption{Iteration count for $\delta = 10^{-5}$}
     \end{subfigure}
     \hspace{1cm}
     \begin{subfigure}[b]{0.4\textwidth}
         \centering
         \includegraphics[width=\textwidth]{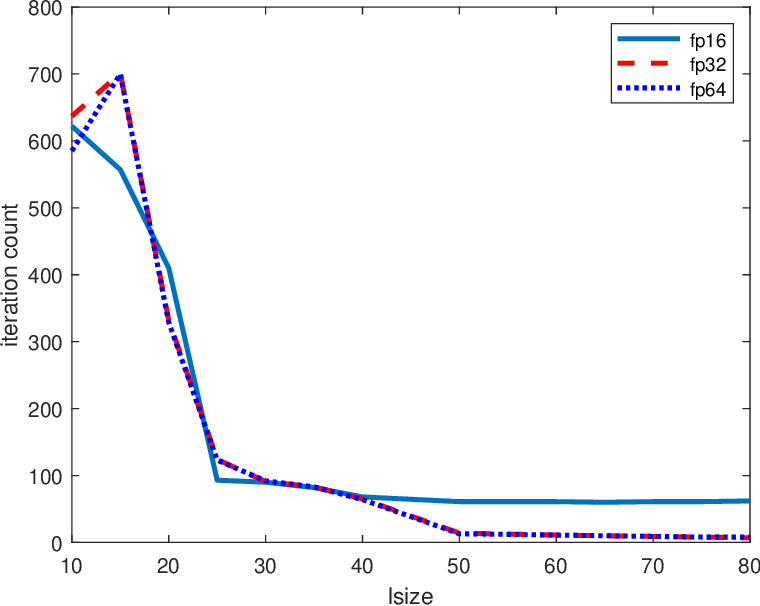}
                \caption{Iteration count for $\delta = 10^{-10}$}
     \end{subfigure}
     
        \caption{The shift $\alpha$ needed to prevent breakdown (top) and the LSQR iteration counts (bottom) for problem mod2 as the {\tt HSL\_MI35} parameter {\tt lsize} that controls the number of entries in the IC factor is increased. The IC factorization is computed using fp16, fp32 and fp64 arithmetic.}
        \label{fig:mod2 alpha}
\end{figure}

Increasing {\tt lsize} can improve the preconditioner
quality because of the resulting reduction in the number of B1 breakdowns and
hence the size of the shift $\alpha$ that is needed for a breakdown-free factorization. Breakdown is more likely when using low precision and a larger $\alpha$ is typically needed.
This is illustrated in Figure~\ref{fig:mod2 alpha} for
problem mod2. We see that for ${\tt lsize}>40$, there are
no breakdowns if fp32 is used and increasing {\tt lsize} beyond 50 leads to no further reductions in the LSQR iteration count.
However, when using fp16 arithmetic, we always get breakdown because the low precision scaled normal matrix $B_l$ is not positive definite, leading to a nonzero shift and higher iteration counts for the fp16 preconditioner compared to the fp32 one. Note that, if we have a fixed amount memory available for $\widetilde L$ and we use fp32 then we can choose a larger {\tt lsize} value than for fp64. This may
result in a lower iteration count for the fp32 version (and hence a saving in the
number of possibly expensive applications of $A$ and $A^T$).

\begin{figure}

          \begin{subfigure}[b]{0.4\textwidth}
         \centering
         \includegraphics[width=\textwidth]{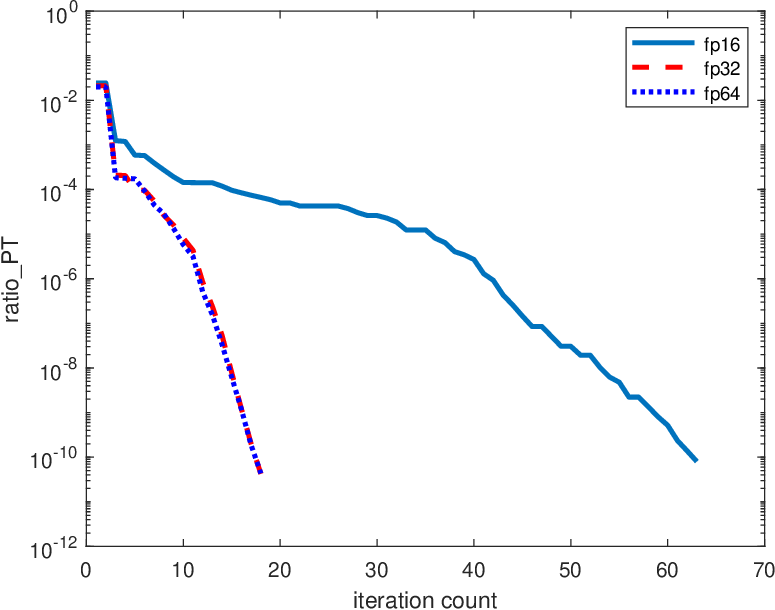}
       \caption{Problem stat96v2}
     \end{subfigure}
     \hspace{1cm}
     \begin{subfigure}[b]{0.4\textwidth}
         \centering
         \includegraphics[width=\textwidth]{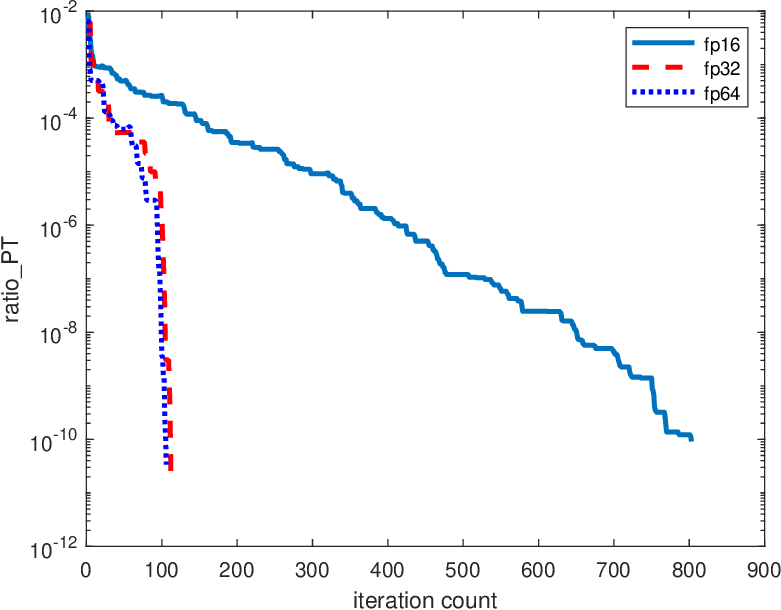}
                \caption{Problem rail2586}
     \end{subfigure}
     
        \caption{Convergence curves for preconditioned LSQR for problems stat96v2 and rail2586. The preconditioner is computed using {\tt HSL\_MI35}
        with the parameter {\tt lsize} set to 20 and 60, respectively (see Table~\ref{T:mi35 fpvar}). $ratio_{PT}$ given by (\ref{eq:ratio_PT}) is plotted against the number of LSQR iterations.}
        \label{fig:stat96v2}
\end{figure}

In Figure~\ref{fig:stat96v2}, LSQR convergence curves
are plotted for problems stat96v2 and rail2586. It is clear that the performance of the preconditioner
computed using fp32 and fp64 arithmetic is comparable, while the curve when using the fp16 preconditioner is significantly delayed. 

So far, we have computed the matrix-vector products with $A$
and $A^T$
and applied the preconditioner using fp64 arithmetic.
For the latter, we  cast the low precision factor entries
to double precision on-the-fly and carry out the operations in double precision.
In Table~\ref{T:mi35 three precisions}, we  use lower
precision for these operations (but still compute the stopping criteria using fp64 arithmetic). Note that, applying the preconditioner in low
precision will produce rounding errors (of order the low precision unit roundoff) that are different at each iteration, so that the application of the
preconditioner is non-constant.
Comparing the results in
Table~\ref{T:mi35 fpvar} with those in Table~\ref{T:mi35 three precisions}, we see that we can successfully use fp32 for computing and applying the preconditioner as well as
for products with $A$  and $A^T$ (that is, we can use fp32 arithmetic  throughout).
If the factors are computed using fp16 then we
can apply the preconditioner using fp16 and then use fp32 for the products with $A$ and $A^T$. With the stopping tolerance $\delta = 10^{-5}$
this can be successful but more iterations
may be required compared to applying the preconditioner
using fp64. A more serious issue is that some problems (pilot\_ja and pilotnov) suffer breakdown, that is, an application of the preconditioner
using fp16 arithmetic results in overflow and the computation
terminates. Thus, although using fp16 
for applying the preconditioner has the potential 
to reduce the cost, in practice checks must be made for
possible breakdowns and, if this happens, entries must be perturbed or a switch made to
applying the preconditioner in higher precision to
ensure the software is robust. We again note that, in some instances, using reorthogonalization within  LSQR can significantly reduce the iteration count.

\begin{center}
\begin{table}[htbp]
{
\footnotesize
\begin{tabular}{lr rl |rrrl}
\hline\Tstrut
 & &
\multicolumn{2}{c|}{fp16/fp16/fp32} &
\multicolumn{4}{c}{fp32/fp32/fp32} \\

\multicolumn{1}{c}{\Tstrut Identifier} &
\multicolumn{1}{c}{{\tt lsize}} &
\multicolumn{2}{c|}{$10^{-5}$} &
\multicolumn{1}{c}{$10^{-5}$} &
\multicolumn{1}{c}{$10^{-10}$} &
\multicolumn{2}{c}{$10^{-15}$}\\

\hline\Tstrut

co9       &  45 &   63 & (39) & 15~~~ & 95~~~ &160 &(141) \\
d2q06c    &  60 &   8 & (8) & 3~~~ & 10~~~ & 15 & (14) \\
delf000   &  20 &  14 & (14) & 4~~~ & 17~~~ & 26 & (25) \\
GE        &  15 &  4 &  (4) & 17~~~ & 27~~~ & 53 & (50) \\
IG5-15    &  60 &  100 &  (96) & 99~~~ & 311~~~ & 476 & (429) \\
illc1033  &  10 &  951 &  (152) &15~~~ & 26~~~ & 38 & (23) \\
illc1850  &  10 &   185 & (80) & 31~~~ & 61~~~ & 89 & (41) \\
Kemelmacher& 30 &  38 &  (38) & 43~~~ & 84~~~ & 132 & (76) \\
large001  &  15 &  27 &  (27) & 4~~~ & 18~~~ & 30 & (26)  \\
mod2      &  30 &  7 &  (10) & 7~~~ & 94~~~ & 165 & (148) \\
pilot\_ja  &  60 &  $\ddag$ & $\ddag$  & 2~~~ & 5~~~ & 10 & (10)  \\
pilotnov  &  60 &  $\ddag$ & $\ddag$   & 2~~~ & 3~~~ & 7 & (7)  \\
psse0     &   5 &  21 & (21) & 5~~~  & 109~~~ & 202 & (153) \\ 
psse1     &  60 &  9 & (9) & 5~~~ & 64~~~ & 116 & (116) \\ 
rail2586  &  60 & 381 & (288) & 109~~~ & 186~~~ & 323 & (95) \\
stat96v2  &  20 & 37 & (37) & 129~~~ & 18~~~ & 29 & (24) \\
watson\_1 &  15 & 79 & ~$\ast$ & 10~~~ & 146~~~ & 235 & ~$\ast$ \\
well1033  &  10 &  12 & (12) & 3~~~ & 5~~~ & 6 & (6) \\
well1850  &  10 &  12 & (11) & 11~~~ & 21~~~ & 27 & (24) \\
world     &  50 &  8 & (8) & 3~~~ & 17~~~ & 26 & (25) \\
\hline
\end{tabular}
}
\caption{Iteration counts for LSQR preconditioned by the
IC factor $L$ computed using {\tt HSL\_MI35}.
fp1/fp2/fp3 denotes $L$ is computed using fp1 arithmetic,
applications of the preconditioner are performed using fp2 arithmetic
and matrix-vector products with $A$ and $A^T$ are performed using fp3 arithmetic.
The stopping criteria (\ref{eq:ratio_PT}) is used
with $\delta = 10^{-5}, 10^{-10}, 10^{-15}$ and is satisfied in each experiment. The numbers in parentheses are the iteration
counts for LSQR with full 
one-sided reorthogonalization. $\ast$ indicates 
insufficient memory and $\ddag$ denotes breakdown. 
}
\label{T:mi35 three precisions}
\end{table}
\end{center}

\subsection{Iterative refinement variants}
Having reported on using LSQR with low precision IC preconditioners, in
Table~\ref{T:lsqr-ir mi35} we compare LSQR, LSQR-IR and GMRES-IR with the IC factorization preconditioner computed using {\tt HSL\_MI35} in single precision. 
The initial solution is set to $x^{(1)} = 0$ (Step 3 of Algorithm~\ref{alg:lsqr-ir}).
We see that LSQR typically has the
smallest $nsol$ and, for some problems, LSQR-IR has a significantly higher count. In our tests, GMRES-IR requires only two outer iterations, but each iteration within the GMRES algorithm applied to the correction equation incorporates orthogonalization and so GMRES-IR is more expensive than simply using LSQR.
If we reduce the stopping tolerance for the outer iteration
(that is, we use a smaller $\delta_2$ in (\ref{eq:ratio_GS}))
then for some examples,  $ratio_{GS}$ reduces before stagnating.
For LSQR-IR it typically stagnates in the range $10^{-9}$ to $10^{-11}$, while for
GMRES-IR it stagnates about $10^{-14}$.

\begin{center}
\begin{table}[htbp]
{
\footnotesize
\begin{tabular}{lrrl rrl rrr}
\hline 
\multicolumn{1}{c}{\Tstrut Identifier} &
\multicolumn{3}{c}{LSQR} &
\multicolumn{3}{c}{LSQR-IR} &
\multicolumn{3}{c}{GMRES-IR} \\
&
\multicolumn{1}{c}{$nsol$} &
\multicolumn{1}{c}{$ratio_{PT}$} &
\multicolumn{1}{c}{$ratio_{GS}$} &
\multicolumn{1}{c}{$nout$} &
\multicolumn{1}{c}{$nsol$} &
\multicolumn{1}{c}{$ratio_{GS}$} &
\multicolumn{1}{c}{$nout$} &
\multicolumn{1}{c}{$nsol$} &
\multicolumn{1}{c}{$ratio_{GS}$} \\
\hline\Tstrut
co9  &        93 &   6.759$\times 10^{-11}$ &  3.611$\times 10^{-8}$ & 8 &   136 &   2.486$\times 10^{-8}$ &     2 &   126 &   1.008$\times 10^{-10}$ \\ 
d2q06c   &    10 &   1.509$\times 10^{-11}$ &  1.253$\times 10^{-9}$ & 5 &   14 &   1.355$\times 10^{-9}$ &     2 &    14 &   1.777$\times 10^{-12}$ \\
delf000  &    17 &   8.620$\times 10^{-11}$ &  1.607$\times 10^{-8}$ & 7 &   26 &   6.593$\times 10^{-9}$ &     2 &    21 &   3.989$\times 10^{-11}$ \\ 
GE       &  27 &     7.334$\times 10^{-11}$ &  4.443$\times 10^{-6}$ & 11 &  54 &   3.633$\times 10^{-9}$ &     2 &    42 &   1.163$\times 10^{-10}$ \\ 
IG5-15    &   305 &  9.505$\times 10^{-11}$ &  3.173$\times 10^{-9}$ & 8 &  396 &   4.456$\times 10^{-9}$ &     2 &   319 &   1.390$\times 10^{-10}$ \\ 
illc1033     & 19 &  7.632$\times 10^{-11}$ &  3.374$\times 10^{-12}$ & 2 &  18 &   1.059$\times 10^{-9}$ &     2 &    20 &   1.454$\times 10^{-11}$ \\
illc1850     & 36 &   5.177$\times 10^{-11}$ &  8.785$\times 10^{-11}$ & 4 & 37 &   1.009$\times 10^{-9}$ &     2 &    60 &   1.584$\times 10^{-11}$ \\ 
Kemelmacher & 68 &   2.344$\times 10^{-11}$ &  6.743$\times 10^{-13}$ &  8 & 218 &   7.986$\times 10^{-9}$ &     2 &    64 &   8.523$\times 10^{-11}$ \\ 
large001 &    17 &   5.106$\times 10^{-11}$ &  1.253$\times 10^{-8}$ &   4 &  18 &   4.746$\times 10^{-9}$ &     2 &    20 &   5.847$\times 10^{-11}$ \\ 
mod2  &       89 &   9.631$\times 10^{-11}$ &  6.997$\times 10^{-8}$ & 8   & 132 &   6.398$\times 10^{-9}$ &     2 &   105 &   1.298$\times 10^{-10}$ \\ 
pilot\_ja  &    5 &   1.638$\times 10^{-11}$ &  6.992$\times 10^{-8}$ &  5 &  10 &  4.256$\times 10^{-10}$ &    2 &     9 &   8.717$\times 10^{-12}$ \\
pilotnov &     3 &   3.769$\times 10^{-11}$ &  2.046$\times 10^{-5}$ & 4   &   8 &  7.232$\times 10^{-10}$ &    2 &     5 &   3.939$\times 10^{-11}$ \\ 
psse0 &      107 &   9.446$\times 10^{-11}$ &  1.502$\times 10^{-9}$ &  10 & 151 &  7.029$\times 10^{-9}$ &     2 &   108 &   9.672$\times 10^{-11}$ \\ 
psse1 &       64 &   7.791$\times 10^{-11}$ &  7.687$\times 10^{-9}$ & 11 &  109 &  6.082$\times 10^{-9}$ &     2 &    80 &   1.948$\times 10^{-10}$ \\
rail2586  &  111 &   4.844$\times 10^{-11}$ &  7.396$\times 10^{-11}$ & 5 &  187 &  4.260$\times 10^{-8}$ &     2 &   360 &   9.319$\times 10^{-11}$ \\ 
stat96v2  &    18 &   4.314$\times 10^{-11}$ &  5.108$\times 10^{-12}$ & 4 &  21 &  1.740$\times 10^{-9}$ &     2 &    20 &   3.899$\times 10^{-11}$ \\
watson\_1 &   141 &   8.930$\times 10^{-11}$ &  5.569$\times 10^{-10}$ & 7 &  148 & 6.938$\times 10^{-9}$ &     2 &   148 &   1.213$\times 10^{-10}$ \\
well1033  &     5 &   1.115$\times 10^{-13}$ &  1.886$\times 10^{-15}$ & 1 &    3 & 1.424$\times 10^{-13}$ &    2 &     3 &   8.953$\times 10^{-14}$ \\
well1850  &    18 &   7.205$\times 10^{-11}$ &  2.882$\times 10^{-11}$ & 4 &   17 & 9.329$\times 10^{-9}$ &     2 &    22 &   2.619$\times 10^{-11}$ \\ 
world     &    17 &   4.657$\times 10^{-11}$ &  1.704$\times 10^{-8}$ & 6 &    24 & 2.229$\times 10^{-9}$ &     2 &    22 &   1.034$\times 10^{-10}$ \\

\hline
\end{tabular}
}
\caption{A comparison of LSQR, LSQR-IR and GMRES-IR. Here {\tt HSL\_MI35} in fp32 arithmetic
is used to compute the IC factorization preconditioner. 
Products with $A$ and $A^T$ and applications of the preconditioner 
are performed using fp64.
For LSQR, (\ref{eq:ratio_PT}) is used with $\delta = 10^{-10}$.
For LSQR-IR, for the correction
equation  (\ref{eq:ratio_PT}) is used with $\delta = 10^{-5}$. 
The convergence tolerance for GMRES applied to the normal
equations for the correction equation is also $10^{-5}$. 
The outer iteration of LSQR-IR and GMRES-IR is terminated using (\ref{eq:ratio_GS})
with $\delta_2  = 10^{-8}$ and (\ref{eq:stop outer})
with $\eta = 10 \times u_{64}$.
$nsol$ denotes the number of solves with $L$ and $L^T$
and $nout$ is the number of outer iterations. $ratio_{GS}$ and $ratio_{PT}$
are given by (\ref{eq:ratio_GS}) and (\ref{eq:ratio_PT}), respectively. }
\label{T:lsqr-ir mi35}
\end{table}
\end{center}

Our next experiment  uses LSQR-IR with the matrix-vector products with $A$ and $A^T$ and applications of the IC factorization preconditioner 
 performed using fp32. Again, the preconditioner is computed using {\tt HSL\_MI35} in fp32 arithmetic. For the correction
equation, the backward error-based stopping condition
(\ref{eq:ratio_PT}) is used with $\delta = 10^{-5}$ and the outer iteration is terminated using (\ref{eq:ratio_GS})
with $\delta_2  = 10^{-8}$. We also impose a limit $itmax=30$ on the number of refinement iterations.
We ran each of our test problems but report only
on a subset in Table~\ref{T:lsqr-ir mi35 stagnation} as these are illustrative of all the results. We observe that, without the test (\ref{eq:stop outer}) for stagnation,
the requested accuracy was not achieved (columns 2 to 4). We therefore use 
(\ref{eq:stop outer})  with $\eta = 10 \times u_{64}$ to monitor
for stagnation. If stagnation occurs, we switch to computing the matrix-vector products and applying the preconditioner 
using fp64; see columns 5 to 9, with the figures in parentheses  reporting the
number of solves with $L$ and $L^T$ in fp64 and number of outer iterations
after the switch to using fp64. Results are included for
using fp64 for all the matrix products and applications of the preconditioner (columns 10 to 12). We see that most of the computation can be performed using fp32 and
then switching to fp64 allows the requested accuracy to be achieved using a relatively
small number of  further iterations.

\begin{center}
\begin{table}[htbp]
{
\footnotesize
\begin{tabular}{lrrl|rlrll|rrl}
\hline 
\multicolumn{1}{c}{\Tstrut Identifier} &
\multicolumn{3}{c|}{fp32/fp32} &
\multicolumn{5}{c|}{fp32/fp32 then fp32/fp64} &
\multicolumn{3}{c}{fp32/fp64} \\
&
\multicolumn{1}{c}{$nout$} &
\multicolumn{1}{c}{$nsol$} &
\multicolumn{1}{c|}{$ratio_{GS}$} &
\multicolumn{2}{c}{$nout$} &
\multicolumn{2}{c}{$nsol$} &
\multicolumn{1}{c|}{$ratio_{GS}$}  &
\multicolumn{1}{c}{$nout$} &
\multicolumn{1}{c}{$nsol$} &
\multicolumn{1}{c}{$ratio_{GS}$} \\
\hline\Tstrut
co9             &    NC &   198 &   2.820$\times 10^{-7}$ & 10 &     (4) &   134 &    (36) &   6.880$\times 10^{-9}$ & 8 &   136 &   2.486$\times 10^{-8}$ \\ 
IG5-15          &    NC &   457 &   1.233$\times 10^{-7}$ &  9 &     (3) &   405 &    (68) &   4.285$\times 10^{-9}$ & 8 &  396 &   4.456$\times 10^{-9}$\\ 
mod2            &    NC &   188 &   2.290$\times 10^{-8}$ &  8 &     (1) &   138 &    (17) &   6.344$\times 10^{-9}$ & 8   & 132 &   6.398$\times 10^{-9}$ \\ 
psse0           &    NC &   242 &   4.082$\times 10^{-8}$ & 12 &     (2) &   187 &    (39) &   1.471$\times 10^{-9}$ &  10 & 151 &  7.029$\times 10^{-9}$ \\

\hline
\end{tabular}
}
\caption{Example results for LSQR-IR. Here {\tt HSL\_MI35} in fp32 arithmetic
is used to compute the IC factorization preconditioner. 
The stopping condition for the correction
equation is (\ref{eq:ratio_PT}) with $\delta = 10^{-5}$. 
The outer iteration is terminated using (\ref{eq:ratio_GS})
with $\delta_2  = 10^{-8}$ and, for the results in columns 5-9, the
stagnation test (\ref{eq:stop outer})
with $\eta = 10 \times u_{64}$. For the results in columns 2-4, fp32 is used for
matrix-vector products and applying the preconditioner. 
For the results in columns 5-9, fp32 is used until stagnation occurs; fp64
is then used for matrix-vector products and applying the preconditioner.
$nsol$ denotes the number of solves with $L$ and $L^T$.
$nout$ is the number of outer iterations. 
In columns 6 and 8, the numbers in parentheses give the
number of solves and outer iterations performed after the
switch to fp64. NC denotes the outer iterations did not converge
within the limit of 30 iterations. $ratio_{GS}$ 
is given by (\ref{eq:ratio_GS}).}
\label{T:lsqr-ir mi35 stagnation}
\end{table}
\end{center}

A possible motivation for using LSQR-IR or GMRES-IR is that it involves 
solving the correction equations with a larger tolerance than is used by LSQR without refinement. We have already seen that 
if we run LSQR
with an fp16 preconditioner and a stopping tolerance of  $10^{-5}$
then (with a random vector $b$) we typically obtain 
the requested accuracy in a relatively small number of iterations
(column 2 of Table~\ref{T:mi35 three precisions}).
This suggests we might expect to use the fp16 preconditioner to successfully 
solve the correction equations within LSQR-IR using only a few iterations,
leading to a small total iteration count for LSQR-IR. However, in practice this is not observed.
When solving the correction equations with LSQR or GMRES, 
the iteration count can be high (significantly higher than
the count needed to get the initial solution).
It is well known that the convergence of Krylov subspace methods  depends
strongly on the right hand side vector. For the correction equation, this vector is $A^T r^{(i)}$. Results are given in Table~\ref{T:lsqr-ir psse0} for problem psse0
with the IC factorization computed using fp16 and applied using fp64.
\begin{center}
\begin{table}[htbp]
{\footnotesize
\begin{tabular}{l rrrrrr}
\hline\Tstrut
Outer  iteration &
1 & 2 & 3 & 4 & 5 & 6 \\
$nsol$ & 16 &      232   &   202  &   333  &  688  &  1254    \\ 
$ratio_{GS}$ &  7.162$\times 10^{-3}$  &   7.030$\times 10^{-4}$  &   1.786$\times 10^{-4}$  &   9.648$\times 10^{-5}$  &   2.747$\times 10^{-5}$ &    5.896$\times 10^{-6}$  \\ 
 \hline\Tstrut
 Outer  iteration &

7 & 8 & 9 & 10 & 11& 12 \\ 
$nsol$ & 387 &   1208 &     393 &    1304   &  66  &  267  \\ 
$ratio_{GS}$ &  1.871$\times 10^{-6}$ &  5.733$\times 10^{-7}$ &  2.611$\times 10^{-7}$ &     6.123$\times 10^{-8}$ & 1.760$\times 10^{-8}$  &   7.298$\times 10^{-9}$ \\ 
\hline
\end{tabular}
}
\caption{Results for LSQR-IR preconditioned by {\tt HSL\_MI35} run in fp16 arithmetic on test problem psse0.
For each refinement iteration, we report the number $nsol$ of
LSQR iterations and $ratio_{GS}$ for the corrected solution. The LSQR stopping tolerance for each correction equation is $10^{-5}$.}
\label{T:lsqr-ir psse0}
\end{table}
\end{center}

\section{Concluding remarks}
\label{sec:conclusions}
In this paper, we have explored the potential for using low precision incomplete
factorization preconditioners for solving linear least-squares
problems. Such preconditioners are important in part because 
there remains a lack of robust general-purpose preconditioners
for solving tough LS problems. We have focused on two approaches: level-based $IC(\ell)$ factorizations and memory-limited factorizations. Our experiments have
been carried out using Fortran code implemented using fp16, fp32 and fp64 arithmetic.
The main findings are summarised as follows.
\begin{itemize}
\item We have demonstrated that the ratio
(\ref{eq:ratio_PT}) based on the estimated quantities as in \cite{pati:2024} provides an effective stopping criteria for
terminating LSQR, which can be used when a 
preconditioner is employed. In the future, we
plan to develop mixed precision implementations of LSQR and CGLS  that
offer the option of employing this stopping criteria for
inclusion in the HSL Library. The software
will be written in Fortran and interfaces to other languages will be provided (including Python and Matlab).
    \item $IC(\ell)$ preconditioners are known to be useful for SPD systems arising from finite difference stencils but our empirical experience is that they are not effective for LS problems, even if the computation is carried out using double precision throughout.
    \item Memory-limited preconditioners can successfully be used
    to solve highly ill-conditioned problems. Provided steps are taken to avoid breakdown, 
    they can be computed in fp16
    arithmetic. In this case, the LSQR iteration counts may be much higher compared to fp32 or fp64 arithmetic (although for some problems, incorporating reorthogonalization can help mitigate this). The higher counts are partly a result of needing
    to employ a larger shift when using fp16 to avoid breakdown. However, if we do not require high accuracy in the LS solution (or we can only afford to perform a few iterations) and/or we only allow a small number of entries in the IC factor $\tilde L$ (for instance, if memory constraints restrict the number of entries or highly sparse factors are sought to reduce the cost of applying the preconditioner) then, taking into account the memory savings, using fp16 is potentially attractive. Furthermore, 
    for LSQR, we have found that using fp32 for the full computation (that is, for computing and applying $\tilde L$ and for products with $A$ and $A^T$) typically performs as well as using fp64 in terms of
    iteration counts and accuracy, while offering memory savings.
\item In recent years there has been significant interest in combining the use of mixed precision with iterative refinement techniques. For least-squares problems using the normal equations we have not been able to demonstrate a consistent advantage in using LSQR-IR or GMRES-IR in place of preconditioned LSQR. However, when the matrix–vector products and preconditioner applications are performed in low precision, LSQR-IR permits a switch to fp64 once stagnation is detected, and this restores the requested accuracy, as illustrated in Table~\ref{T:lsqr-ir mi35 stagnation}.
\end{itemize}

In this study, we do not report timings. As already noted, it is currently
not possible to obtain meaningful timings for our Fortran implementations with the NAG compiler but this is clearly of future interest.
In our test environment, we did not observe speedup using fp32. This is because our problems are relatively small and our implementations are serial. This means that the matrix-vector products with $A$ and $A^T$ and triangular solves are likely to be latency-bound. For large problems (which do not fit into the cache) and a multithreaded implementation (which saturates the bandwidth), the matrix-vector products will be memory-bound and result in a speedup from fp32, particularly on modern architectures for which fp32 is significantly faster than fp64.

Finally, we observe that another version of 16-bit arithmetic, usually referred to as  bfloat16 or bf16, was developed by Google specifically
for deep learning training on their Tensor Product Units.
Currently, no mainstream Fortran compiler supports bfloat16
on CPUs\footnote{NVIDIA's NVFortran, designed for their GPUs, does support bfloat16.}. Should one become available, it would be of interest to compare its performance to that of fp16 when solving LS problems. The attraction is that
 bfloat16 has the same exponent size as fp32 making converting from fp32 to bfloat16 straightforward (overflow and underflow do not occur in the conversion).
The key disadvantage of bfloat16 is its lesser precision: essentially three significant decimal digits versus four for fp16, which may make it unsuitable for solving challenging LS problems.

\bigskip
\noindent {\bf Acknowledgments}\\
We are grateful to two reviewers for the constructive and helpful suggestions that we have incorporated into paper.

\def\cprime{$'$}


\begin{thebibliography}{46}
\providecommand{\natexlab}[1]{#1}
\providecommand{\url}[1]{{#1}}
\providecommand{\urlprefix}{URL }
\providecommand{\doi}[1]{\url{https://doi.org/#1}}
\providecommand{\eprint}[2][]{\url{#2}}
 \bibcommenthead

\bibitem[{Abdelfattah et~al(2021)Abdelfattah, Anzt, Boman, Carson, Cojean,
  Dongarra, Fox, Gates, Higham, Li, Liu, Loe, Luszczek, Pranesh, Rajamanickam,
  Ribizel, Smith, Swirydowicz, Thomas, Tomov, Tzai, Yamazaki, and
  Yang}]{abdel:2021}
Abdelfattah A, Anzt H, Boman E, et~al (2021) A survey of numerical linear
  algebra methods utilizing mixed precision. Int J High Performance Comput Appl
  35(4):344--369

\bibitem[{Amestoy et~al(2024)Amestoy, Buttari, Higham, L'Excellent, Mary, and
  Vieubl{\'e}}]{abhltv:2024}
Amestoy P, Buttari A, Higham NJ, et~al (2024) Five-precision {GMRES}-based
  iterative refinement. SIAM J Matrix Anal Appl 45(1):529--552.
  \doi{10.1137/23M1549079}

\bibitem[{Avron et~al(2019)Avron, Druinsky, and Toledo}]{avdt:2019}
Avron H, Druinsky A, Toledo S (2019) Spectral condition-number estimation of
  large sparse matrices. Numer Linear Algebra Appl 26(3):e2235.
  \doi{10.1002/nla.2235}

\bibitem[{Bj{\"o}rck(1967)}]{bjor:67a}
Bj{\"o}rck {\AA} (1967) Iterative refinement of linear least squares solutions
  {I}. BIT 7(4):257--278. \doi{10.1007/bf01939321}

\bibitem[{Bj{\"o}rck(2024)}]{bjor:2024}
Bj{\"o}rck {\AA} (2024) Numerical Methods for Least Squares Problems, 2nd edn.
  SIAM, Philadelphia, PA, \doi{10.1137/1.9781611977950}

\bibitem[{Bru et~al(2014)Bru, Mar\'{\i}n, Mas, and T\r{u}ma}]{bmmt:2014}
Bru R, Mar\'{\i}n J, Mas J, et~al (2014) Preconditioned iterative methods for
  solving linear least squares problems. SIAM J Sci Comput 36(4):A2002--A2022.
  \doi{10.1137/130931588}

\bibitem[{Carson and Dau\v~zickait\.e(2025{\natexlab{a}})}]{cada:2024}
Carson E, Dau\v~zickait\.e I (2025{\natexlab{a}}) A comparison of mixed
  precision iterative refinement approaches for least-squares problems. SIAM J
  Matrix Anal Appl 46(2):1117--1144. \doi{10.1137/24M1664927},
  \urlprefix\url{https://doi.org/10.1137/24M1664927}

\bibitem[{Carson and Dau\v~zickait\.e(2025{\natexlab{b}})}]{cada:2024a}
Carson E, Dau\v~zickait\.e I (2025{\natexlab{b}}) Mixed precision sketching for
  least-squares problems and its application in {GMRES}-based iterative
  refinement. SIAM J Matrix Anal Appl 46(3):2041--2060.
  \doi{10.1137/24M1702246}, \urlprefix\url{https://doi.org/10.1137/24M1702246}

\bibitem[{Carson and Higham(2017)}]{cahi:17}
Carson E, Higham NJ (2017) A new analysis of iterative refinement and its
  application to accurate solution of ill-conditioned sparse linear systems.
  SIAM J Sci Comput 39(6):A2834--A2856. \doi{10.1137/17M1122918}

\bibitem[{Carson and Higham(2018)}]{cahi:18}
Carson E, Higham NJ (2018) Accelerating the solution of linear systems by
  iterative refinement in three precisions. SIAM J Sci Comput 40(2):A817--A847.
  \doi{10.1137/17M1140819}

\bibitem[{Carson et~al(2020)Carson, Higham, and Pranesh}]{cahp:2000}
Carson E, Higham NJ, Pranesh S (2020) Three-precision {GMRES}-based iterative
  refinement for least squares problems. SIAM J Sci Comput 42(6):A4063--A4083.
  \doi{10.1137/20M1316822}, \urlprefix\url{https://doi.org/10.1137/20M1316822}

\bibitem[{Chang et~al(2009)Chang, Paige, and Titley-Peloquin}]{cptp:2009}
Chang XW, Paige CC, Titley-Peloquin D (2009) Stopping criteria for the
  iterative solution of linear least squares problems. SIAM J Matrix Anal Appl
  31(2):831--852. \doi{10.1137/080724071}

\bibitem[{Fong and Saunders(2011)}]{fosa:2011}
Fong DCL, Saunders M (2011) {LSMR}: An iterative algorithm for sparse
  least-squares problems. SIAM J Sci Comput 33(5):2950--2971.
  \doi{10.1137/10079687X}

\bibitem[{Georgiou et~al(2023)Georgiou, Boutsikas, Drineas, and
  Anzt}]{gbda:2023}
Georgiou V, Boutsikas C, Drineas P, et~al (2023) A mixed precision randomized
  preconditioner for the {LSQR} solver on {GPU}s. In: International Conference
  on High Performance Computing, Springer, pp 164--181

\bibitem[{Gould and Scott(2017)}]{gosc:2017}
Gould NIM, Scott JA (2017) The state-of-the-art of preconditioners for sparse
  linear least-squares problems. ACM Trans Math Software 43(4):Art. 36, 1--35.
  \doi{10.1145/3014057}

\bibitem[{Hallman(2020)}]{hall:2020}
Hallman E (2020) Estimating the backward error for the least-squares problem
  with multiple right-hand sides. Linear Algebra Appl 605:227--238.
  \doi{10.1016/j.laa.2020.07.020}

\bibitem[{Hestenes and Stiefel(1952)}]{hest:52}
Hestenes MR, Stiefel E (1952) Methods of conjugate gradients for solving linear
  systems. J Research of the National Bureau of Standards 49(6):409--436

\bibitem[{Higham(2002)}]{high:02}
Higham NJ (2002) Accuracy and Stability of Numerical Algorithms, 2nd edn. SIAM,
  Philadelphia, PA, \doi{10.1137/1.9780898718027},
  \urlprefix\url{https://doi.org/10.1137/1.9780898718027}

\bibitem[{Higham and Mary(2022)}]{hima:2022}
Higham NJ, Mary T (2022) Mixed precision algorithms in numerical linear
  algebra. Acta Numer 31:347--414. \doi{10.1017/S0962492922000022}

\bibitem[{Higham and Pranesh(2021)}]{hipr:2021}
Higham NJ, Pranesh S (2021) Exploiting lower precision arithmetic in solving
  symmetric positive definite linear systems and least squares problems. SIAM J
  Sci Comput 43(1):A258--A277. \doi{10.1137/19M1298263}

\bibitem[{Higham and Stewart(1987)}]{hist:87}
Higham NJ, Stewart GW (1987) Numerical linear algebra in statistical computing.
  In: Iserles I, Powell MJD (eds) The State of the Art in Numerical Analysis.
  Oxford University Press

\bibitem[{Hogg et~al(2010)Hogg, Reid, and Scott}]{hors:2010}
Hogg JD, Reid JK, Scott JA (2010) Design of a multicore sparse {C}holesky
  factorization using {DAG}s. SIAM J Sci Comput 32(6):3627--3649.
  \doi{10.1137/090757216}, \urlprefix\url{https://doi.org/10.1137/090757216}

\bibitem[{HSL(accessed 2025)}]{hsl:2025}
HSL (accessed 2025) {A} collection of {Fortran} codes for large-scale
  scientific computation. \url{http://www.hsl.rl.ac.uk}

\bibitem[{Hysom and Pothen(2002)}]{hypo:02}
Hysom D, Pothen A (2002) Level-based incomplete {LU} factorization: Graph model
  and algorithms. Preprint UCRL-JC-150789, US Department of Energy

\bibitem[{Jir{\'a}nek and Titley-Peloquin(2010)}]{jitp:2010}
Jir{\'a}nek P, Titley-Peloquin D (2010) Estimating the backward error in
  {LSQR}. SIAM J Matrix Anal Appl 31(4):2055--2074. \doi{10.1137/090770655}

\bibitem[{Kashi et~al(2024)Kashi, Lu, Brewer, Rogers, Matheson, Shankar, and
  Wang}]{klbr:2024}
Kashi A, Lu H, Brewer W, et~al (2024) Mixed-precision numerics in scientific
  applications: survey and perspectives. arXiv preprint arXiv:241219322

\bibitem[{Klein and Lu(1996)}]{kllu:1996}
Klein P, Lu HI (1996) Efficient approximation algorithms for semidefinite
  programs arising from max cut and coloring. In: Proceedings of the
  Twenty-Eighth Annual ACM Symposium on Theory of Computing. Association for
  Computing Machinery, New York, NY, USA, STOC \'96, p. 338--347,
  \doi{10.1145/237814.237980},
  \urlprefix\url{https://doi.org/10.1145/237814.237980}

\bibitem[{Li(2024)}]{li:2024}
Li H (2024) Double precision is not necessary for {LSQR} for solving discrete
  linear ill-posed problems. J Sci Comput 98(3):55.
  \doi{10.1007/s10915-023-02447-4}

\bibitem[{Lin and Mor{\'e}(1999)}]{limo:99}
Lin CJ, Mor{\'e} JJ (1999) Incomplete {C}holesky factorizations with limited
  memory. SIAM J Sci Comput 21(1):24--45. \doi{10.1137/S1064827597327334}

\bibitem[{Meurant et~al(2021)Meurant, Pape\v{z}, and Tich\'y}]{mept:2021}
Meurant G, Pape\v{z} J, Tich\'y P (2021) Accurate error estimation in {CG}.
  Numer Algorithms 88(3):1337--1359. \doi{10.1007/s11075-021-01078-w},
  \urlprefix\url{https://doi.org/10.1007/s11075-021-01078-w}

\bibitem[{Nachtigal et~al(1992)Nachtigal, Reddy, and Trefethen}]{nart:1992}
Nachtigal NM, Reddy SC, Trefethen LN (1992) How fast are nonsymmetric matrix
  iterations? SIAM J Matrix Anal Appl 13(3):778--795. \doi{10.1137/0613049}

\bibitem[{Paige and Saunders(1982)}]{pasa:82}
Paige CC, Saunders MA (1982) {LSQR}: An algorithm for sparse linear equations
  and sparse least squares. ACM Trans Math Software 8(1):43--71.
  \doi{10.1145/355984.355989}

\bibitem[{Pape\v{z} and Tich\'y(2024)}]{pati:2024}
Pape\v{z} J, Tich\'y P (2024) Estimating error norms in {CG}-like algorithms
  for least-squares and least-norm problems. Numer Algorithms 97(1):1--28.
  \doi{10.1007/s11075-023-01691-x},
  \urlprefix\url{https://doi.org/10.1007/s11075-023-01691-x}

\bibitem[{Scott and T\r{u}ma(2011)}]{sctu:2011}
Scott JA, T\r{u}ma M (2011) The importance of structure in incomplete
  factorization preconditioners. BIT 51(2):385--404.
  \doi{10.1007/s10543-010-0299-8}

\bibitem[{Scott and T\r{u}ma(2014{\natexlab{a}})}]{sctu:2014a}
Scott JA, T\r{u}ma M (2014{\natexlab{a}}) {\tt HSL\_MI28}: an efficient and
  robust limited-memory incomplete {C}holesky factorization code. ACM Trans
  Math Software 40(4):Art. 24, 1--19. \doi{10.1145/2617555}

\bibitem[{Scott and T\r{u}ma(2014{\natexlab{b}})}]{sctu:2014}
Scott JA, T\r{u}ma M (2014{\natexlab{b}}) On positive semidefinite modification
  schemes for incomplete {C}holesky factorization. SIAM J Sci Comput
  36(2):A609--A633. \doi{10.1137/130917582}

\bibitem[{Scott and T\r{u}ma(2017)}]{sctu:2017b}
Scott JA, T\r{u}ma M (2017) Solving mixed sparse-dense linear least-squares
  problems by preconditioned iterative methods. SIAM J Sci Comput
  39(6):A2422--A2437. \doi{10.1137/16M1108339}

\bibitem[{Scott and T\r{u}ma(2019)}]{sctu:2019a}
Scott JA, T\r{u}ma M (2019) Sparse stretching for solving sparse-dense linear
  least-squares problems. SIAM J Sci Comput 41:A1604--A1625.
  \doi{10.1137/18M1181353}

\bibitem[{Scott and T\r{u}ma(2022)}]{sctu:2022}
Scott JA, T\r{u}ma M (2022) A computational study of using black-box {QR}
  solvers for large-scale sparse-dense linear least squares problems. ACM Trans
  Math Software 48(1):Art. 5, 1--24. \doi{10.1145/3494527}

\bibitem[{Scott and T\r{u}ma(2024)}]{sctu:2024}
Scott JA, T\r{u}ma M (2024) Avoiding breakdown in incomplete factorizations in
  low precision arithmetic. ACM Trans Math Software 50(2):Art. 9, 1--25.
  \doi{10.1145/3651155}

\bibitem[{Scott and T\r{u}ma(2025)}]{sctu:2025}
Scott JA, T\r{u}ma M (2025) Developing robust incomplete {C}holesky
  factorizations in half precision arithmetic. Numer Algors pp 1--22.
  \doi{10.1007/s11075-025-02015-x}

\bibitem[{Simon(1984)}]{simo:1984}
Simon HD (1984) Analysis of the symmetric {L}anczos algorithm with
  reorthogonalization methods. Linear Algebra Appl 61:101--131.
  \doi{10.1016/0024-3795(84)90025-9}

\bibitem[{Simon and Zha(2000)}]{sizh:2000}
Simon HD, Zha H (2000) Low-rank matrix approximation using the {L}anczos
  bidiagonalization process with applications. SIAM J Sci Comput
  21(6):2257--2274. \doi{10.1137/S1064827597327309}

\bibitem[{Tismenetsky(1991)}]{tism:91}
Tismenetsky M (1991) A new preconditioning technique for solving large sparse
  linear systems. Linear Algebra Appl 154:331--353.
  \doi{10.1016/0024-3795(91)90383-8}

\bibitem[{Watts-III(1981)}]{watt:81}
Watts-III JW (1981) A conjugate gradient truncated direct method for the
  iterative solution of the reservoir simulation pressure equation. Soc
  Petroleum Engrg J 21:345--353. \doi{10.2118/8252-PA}

\bibitem[{Zhang and Wu(2019)}]{zhwu:2019}
Zhang S, Wu P (2019) High accuracy low precision qr factorization and least
  square solver on {GPU} with tensorcore. arXiv preprint arXiv:191205508

\end{thebibliography}

\end{document}